\documentclass{amsart}

\usepackage{amsmath} 
\usepackage{amssymb}

\newtheorem{theorem}{Theorem}[section] 
\newtheorem{claim}{Claim}[theorem]
\newtheorem{lemma}[theorem]{Lemma} 
\newtheorem{proposition}[theorem]{Proposition} 
\newtheorem{observation}[theorem]{Observation} 
\newtheorem{corollary}[theorem]{Corollary} 

\theoremstyle{definition}
\newtheorem{definition}[theorem]{Definition}

\newtheorem{problem}[theorem]{Problem}

\theoremstyle{remark}
\newtheorem{remark}[theorem]{Remark}

\numberwithin{equation}{section}
\newcommand{\forces}{\Vdash}

\newcommand{\bV}{{\bf V}} 
\newcommand{\lesdot}{\mathrel{\mathord{<}\!\!\raise 
0.8 pt\hbox{$\scriptstyle\circ$}}}

\newcommand{\conc}{{}^\frown\!}
\newcommand{\lh}{{\rm lh}\/}
\newcommand{\rest}{{\restriction}}
\newcommand{\dom}{{\rm dom}}

\newcommand{\rk}{{\rm rk}}
\newcommand{\vtl}{\vartriangleleft}

\newcommand{\Mkkk}{{\bf M}^\lambda_{\lambda,\lambda}}

\newcommand{\cA}{{\mathcal A}}

\newcommand{\bbC}{{\mathbb C}}

\newcommand{\cD}{{\mathcal D}}

\newcommand{\cF}{{\mathcal F}}
\newcommand{\cG}{{\mathcal G}}
\newcommand{\bG}{{\bf G}}

\newcommand{\bbL}{{\mathbb L}}

\newcommand{\bbP}{{\mathbb P}}
\newcommand{\cP}{{\mathcal P}}

\newcommand{\bbQ}{{\mathbb Q}}
\newcommand{\dbQ}{{\name{\mathbb Q}}}

\newcommand{\cT}{{\mathcal T}}
\newcommand{\cU}{{\mathcal U}}
\newcommand{\bU}{{\bf U}}

\newcommand{\cX}{{\mathcal X}}

\newcommand{\cf}{{\rm cf}\/} 
 
\newcommand{\fil}{{\rm fil}\/} 
 
\newcommand{\otp}{{\rm otp}\/} 

\newcommand{\st}{{\bf st}} 
\newcommand{\cHchi}{{\mathcal H}(\chi)}
\newcommand{\vare}{\varepsilon}

\newcommand{\bqz}{{{\mathbb Q}^0_\lambda}}
\newcommand{\bqd}{{{\mathbb Q}^*_\lambda}}
\newcommand{\bqo}{{{\mathbb Q}^0_{\omega_1}}}

\newcommand{\pr}{{\rm pr}}

\newcommand{\Gcc}{{\Game^{\rm cc}_{\varepsilon,\lambda}}}

\newcount\skewfactor
\def\mathunderaccent#1#2 {\let\theaccent#1\skewfactor#2
\mathpalette\putaccentunder}
\def\putaccentunder#1#2{\oalign{$#1#2$\crcr\hidewidth
\vbox to.2ex{\hbox{$#1\skew\skewfactor\theaccent{}$}\vss}\hidewidth}}
\def\name{\mathunderaccent\tilde-3 }

\begin{document}

\title{Generating ultrafilters in a reasonable way}

\author{Andrzej Ros{\l}anowski}
\address{Department of Mathematics\\
 University of Nebraska at Omaha\\
 Omaha, NE 68182-0243, USA}
\email{roslanow@member.ams.org}
\urladdr{http://www.unomaha.edu/logic}

\author{Saharon Shelah}
\address{Einstein Institute of Mathematics\\
Edmond J. Safra Campus, Givat Ram\\
The Hebrew University of Jerusalem\\
Jerusalem, 91904, Israel\\
 and  Department of Mathematics\\
 Rutgers University\\
 New Brunswick, NJ 08854, USA}
\email{shelah@math.huji.ac.il}
\urladdr{http://shelah.logic.at}
\thanks{The first author would like to thank the Hebrew University of
  Jerusalem and the Lady Davis Fellowship Trust for awarding him with
  {\em Sch\"onbrunn Visiting Professorship\/} under which this research was
  carried out. \\
Both authors acknowledge support from the United States-Israel
Binational Science Foundation (Grant no. 2002323). This is publication
889 of the second author.}     

\date{September 2007}

\begin{abstract}
  We continue investigations of {\em reasonable ultrafilters\/} on
  uncountable cardinals defined in Shelah \cite{Sh:830}. We introduce a
  general scheme of generating a filter on $\lambda$ from filters on smaller
  sets and we investigate the combinatorics of objects obtained this way.
\end{abstract}

\maketitle

\section{Introduction}
{\em Reasonable ultrafilters\/} were introduced in Shelah \cite{Sh:830} in
order to suggest a line of research that would repeat in some sense the
beautiful theory created around the notion of {\em P--points on
  $\omega$}. If we are interested in generalizing $P$--points, but we do not
want to deal with large cardinals, we have to be somewhat creative in
re-interpreting the property that {\em any countable family of members of
  the ultrafilter has a pseudo-intersection in the ultrafilter}. An
interesting way of doing this is to look at the ways an ultrafilter on an
uncountable cardinal $\lambda$ can be obtained from $\lambda$--sequences of
objects on smaller cardinals. The general scheme for this approach is
motivated by \cite[\S 5,6]{RoSh:470} and it is presented in Definition
\ref{genfil}. In this context the $P$--pointness of an ultrafilter may be
re-interpreted as $(<\lambda^+)$--directness of its generating system (see
\ref{remres}).   

As in \cite{Sh:830}, when working with ultrafilters on $\lambda$, we want to
concentrate on those which are {\em very non-normal}.  Thus very often we
ask ourselves questions concerning {\em weak reasonability\/} of the
ultrafilter obtained from a generating system, and the following property is 
always of interest in this paper.

\begin{definition}
[Shelah {\cite[Def. 1.4]{Sh:830}}]
\label{1.5}
\begin{enumerate}
\item We say that a uniform ultrafilter $D$ on $\lambda$ is {\em weakly
  reasonable\/} if for every non-decreasing unbounded function $f\in
{}^\lambda\lambda$ there is a club $C$ of $\lambda$ such that    
\[\bigcup\{[\delta,\delta + f(\delta)):\delta \in C\}\notin D.\]
\item Let $D$ be an ultrafilter on $\lambda$, $C\subseteq\lambda$ be  a club
  and let $\langle\delta_\xi:\xi<\lambda\rangle$ be the increasing
  enumeration of $C\cup\{0\}$. We define 
\[D/C=\big\{A\subseteq\lambda: \bigcup_{\xi\in A}[\delta_\xi,
\delta_{\xi+1}) \in D\big\}.\]
(It is an ultrafilter on $\lambda$.) $D/C$ will be called {\em the quotient
  of $D$ by $C$}. 
\end{enumerate}
\end{definition}

\begin{observation}
[Shelah {\cite[Obs. 1.5]{Sh:830}}]
\label{easyob}
Let $D$ be a uniform ultrafilter on a regular uncountable cardinal
$\lambda$. Then the following conditions are equivalent: 
\begin{enumerate}
\item[(A)] $D$ is weakly reasonable,
\item[(B)] for every increasing continuous sequence $\langle\delta_\xi:\xi< 
\lambda\rangle\subseteq\lambda$ there is a club $C^*$ of $\lambda$ such that 
\[\bigcup\big\{[\delta_\xi,\delta_{\xi+1}):\xi\in C^*\big\}\notin D,\] 
\item[(C)] for every club $C$ of $\lambda$ the quotient $D/C$ does not
  extend the filter generated by clubs of $\lambda$.
\end{enumerate}
\end{observation}
\medskip

This paper continues Shelah \cite{Sh:830} and Ros{\l}anowski and Shelah
\cite{RoSh:890}, but it is essentially self contained. In the first section we
present our key definitions introducing {\em systems of local filters\/} and
partial orders $\bqd(\cF)$, $\bqz(\cF)$ associated with them. We explain how
those partial orders can be made $(<\lambda^+)$--complete (in
\ref{whenclosed}, \ref{closureclosed}) and we show that ultrafilters
generated by sufficiently directed generating systems are weakly reasonable,
unless they are produced from a measurable ultrafilter (see
\ref{nonres}). The second section is concerned with the {\em full system
$\cF^{\rm ult}$ of local ultrafilters\/} and the ultrafilters on $\lambda$
generated by $H\subseteq\bqd(\cF^{\rm ult})$. We show that there may be
weakly reasonable ultrafilters on $\lambda$ generated by some
$H'\subseteq\bqz(\cF)$ which cannot be obtained by use of $\cF^{\rm ult}$
(see \ref{nonultra}). Furthermore, we introduce more properties of families
$H \subseteq\bqd(\cF^{\rm ult})$ which are useful in generating
ultrafilters on $\lambda$. In the third section we are interested in a
system $\cF^{\rm pr}$ of local filters and its relation to generating
numbers (in standard sense) of filters on $\lambda$ (see \ref{iterdom},
\ref{criterion}). Finally, in the last section we show that the
inaccessibility of $\lambda$ in the assumptions of
\cite[Prop. 1.6(1)]{Sh:830} is needed: consistently, there is a very
reasonable ultrafilter $D$ on $\omega_1$ such that Odd has a winning
strategy in $\Game_D$ (see \ref{inacneeded}).  
\medskip

\noindent {\bf Notation:}\quad Our notation is rather standard and
compatible with that of classical textbooks (like Jech \cite{J}). 

\begin{enumerate}
\item Ordinal numbers will be denoted be the lower case initial letters of
the Greek alphabet ($\alpha,\beta,\gamma,\delta\ldots$) and also by $i,j$
(with possible sub- and superscripts). Cardinal numbers will be called
$\kappa,\lambda,\mu$ (with possible sub- and superscripts). {\bf $\lambda$
  is always assumed to be an uncountable regular cardinal}. 

\item For two sequences $\eta,\nu$ we write $\nu\vartriangleleft\eta$
whenever $\nu$ is a proper initial segment of $\eta$, and $\nu
\trianglelefteq\eta$ when either $\nu\vartriangleleft\eta$ or $\nu=\eta$. 
The length of a sequence $\eta$ is denoted by $\lh(\eta)$.

\item We will use letters $D,E,F$ and $d$ (with possible indexes) to denote
  filters on various sets. Typically, $D$ will be a filter on $\lambda$
  (possibly an ultrafilter), while $E,F$ will stand for filters on smaller
  sets. Also, in most cases $d$ will be an ultrafilter on a set of size less
  than $\lambda$.\\
For a filter $F$ of subsets of a set $A$, the family of all $F$--positive
subsets of $A$ is called $F^+$. (So $B\in F^+$ if and only if $B\subseteq A$
and $B\cap C\neq\emptyset$ for all $C\in F$.) 

\item In forcing we keep the older convention that {\em a stronger condition
    is the larger one}.  For a forcing notion $\bbP$, $\Gamma_\bbP$ stands
  for the canonical $\bbP$--name for the generic filter in $\bbP$. With this
  one exception,  all $\bbP$--names for objects in the extension via $\bbP$
  will be denoted with a tilde below (e.g., $\name{\tau}$, $\name{X}$). 
\end{enumerate}

\section{Generating a filter from systems of local filters}
Here we present the general scheme of generating a filter on a {\bf regular
  uncountable} cardinal $\lambda$ by using smaller filters. Our approach is
slightly different from the one in \cite[\S 2]{Sh:830} and/or \cite[\S
1]{RoSh:890}, but the difference is notational only (see \ref{remres} below).

\begin{definition}
  \label{bases}
  \begin{enumerate}
  \item A {\em system of local filters on $\lambda$\/} is a family $\cF$
    such that 
    \begin{itemize}
    \item all members of $\cF$ are triples $(\alpha,Z,F)$ such that
      $Z\subseteq \lambda$, $|Z|<\lambda$, $\alpha=\min(Z)$ and $F$ is a
      proper filter on $Z$, 
     \item the set $\big\{\alpha<\lambda:\big(\exists Z,F\big)
       \big((\alpha,Z,F)\in \cF\big)\big\}$ is unbounded in $\lambda$. 
    \end{itemize}
If above for every $(\alpha,Z,F)\in\cF$, the set $Z$ is infinite and $F$ is
a non-principal ultrafilter on $Z$, then we say that {\em $\cF$ is a system
of  local non-principal ultrafilters}.  
  \item More generally, if $\Psi$ is a property of filters, then a {\em system
      of local $\Psi$--filters on $\lambda$\/} is a system of local filters
    $\cF$ such that for every $(\alpha,Z,F)\in\cF$, the filter $F$ has the
    property $\Psi$.  The {\em full system of local $\Psi$--filters\/} is the
    family of {\em all\/} triples $(\alpha,Z,F)$ such that $\alpha<\lambda$,
    $\alpha\in Z\subseteq \lambda\setminus \alpha$, $|Z|<\lambda$  and $F$
    is a proper filter on $Z$ with the property $\Psi$ (assuming that it
    forms a system of local filters). The full system of local non-principal
    ultrafilters on $\lambda$ is denoted by $\cF^{\rm ult}_\lambda$ or just
    $\cF^{\rm ult}$ (if $\lambda$ is understood).  
  \end{enumerate}
\end{definition}

The next definition introduces the filters generated by some families of
local filters. As we have said in the introduction, our motivations have
roots in forcings with norms and this suggested us to use sometimes a
forcing-like notation (e,g, $\bqd$) similar to that of \cite{RoSh:470}. It is
also worth noticing that some families of generators may be used as forcing
notions - for instance $(\bqz,\leq^*)$ is the forcing used in the end of
\cite[Sec. 1]{RoSh:890}.  

\begin{definition}
\label{genfil} 
Let $\cF$ be a system of local filters on $\lambda$. 
\begin{enumerate}
\item We let $\bqd(\cF)$ be the family of all sets $r\subseteq \cF$  such that
 \[\big(\forall\xi<\lambda\big)\big(|\{(\alpha,Z,F)\in r:\alpha=\xi\}|<
 \lambda\big)\quad\mbox{ and }\quad |r|=\lambda.\]  
For $r\in\bqd(\cF)$ we define
\[\fil(r)=\big\{A\subseteq\lambda:\big(\exists\vare<\lambda\big)\big(\forall 
(\alpha,Z,F)\in r\big)\big(\vare\leq\alpha\ \Rightarrow\ A \cap Z\in F\big)
\big\},\] 
and we define a binary relation $\leq^*=\leq^*_{\cF}$ on $\bqd(\cF)$ by  

$r_1\leq^*_{\cF} r_2$ if and only if ($r_1,r_2\in\bqd(\cF)$ and)
$\fil(r_1)\subseteq \fil(r_2)$.
\item We say that an $r\in\bqd(\cF)$ is {\em strongly disjoint\/} if and
  only if  
\begin{itemize}
\item $\big(\forall\xi<\lambda\big)\big(|\{(\alpha,Z,F)\in r:\alpha=\xi\}|   
<2\big)$, and 
\item $\big(\forall (\alpha_1,Z_1,F_1),(\alpha_2,Z_2,F_2)\in r\big)\big(
\alpha_1<\alpha_2\ \Rightarrow\ Z_1\subseteq\alpha_2\big)$.
\end{itemize}
We let $\bqz(\cF)$ be the collection of all strongly disjoint elements of
$\bqd(\cF)$. 
\item We write $\bqd,\bqz$ for  $\bqd(\cF^{\rm ult}),\bqz(\cF^{\rm ult})$,
  respectively (where, remember, $\cF^{\rm ult}$ is the full system of local 
  non-principal ultrafilters). 
\item For a set $H\subseteq\bqd(\cF)$ we let $\fil(H)=\bigcup
  \big\{\fil(r):r\in H\}$.   
\end{enumerate}
\end{definition}

\begin{remark}
\label{remres}
\begin{enumerate}
\item Note that if $r\in \bqz$ then there is $r'\in\bqz$ such that
$\fil(r')=\fil(r)$ and for some club $C$ of $\lambda$ we have  
\[\big\{(\alpha,Z):\big(\exists d\big)\big((\alpha,Z,d)\in r'\big)\big\}= 
\big\{(\alpha,[\alpha,\beta)):\alpha\in C\ \&\ \beta=\min\big(C\setminus 
(\alpha+1)\big)\big\}.\] 
Thus $\bqz$ is essentially the same as the one defined in
\cite[Def. 2.5]{Sh:830}. 
\item If $H\subseteq\bqd(\cF)$ is $\leq^*$--directed, then $D=\fil(H)$ 
  is a filter on $\lambda$ extending the filter of co-bounded sets. We may
  say the that {\em the filter $D$ is generated by $H$} or that {\em $H$
    is the generating system for $D$.} 
\end{enumerate}
\end{remark}

\begin{definition}
\label{eproduct} 
Suppose that
\begin{enumerate}
\item[(a)] $X$ is a non-empty set and $F$ is a filter on $X$,
\item[(b)] $F_x$ is a filter on a set $Z_x$ (for $x\in X$). 
\end{enumerate}
We let 
\[\bigoplus\limits^F_{x\in X} F_x=\big\{A\subseteq\bigcup\limits_{x\in X}
Z_x: \{x\in X:Z_x\cap A\in F_x\}\in F\big\}.\]
(Clearly, $\bigoplus\limits^F_{x\in X} F_x$ is a filter on
$\bigcup\limits_{x\in X} Z_x$.) If $X$ is a linearly ordered set (e.g. it is
a set of ordinals) with no maximal element and $F$ is the filter of
all co-bounded subsets of $X$, then we will write $\bigoplus\limits_{x\in X} 
F_x$ instead of $\bigoplus\limits^F_{x\in X} F_x$.  
\end{definition}

\begin{proposition}
[Cf. {\cite[Prop. 2.9]{Sh:830}}]
\label{1.3C}
\begin{enumerate}
\item Let $\cF$ be a system of local filters on $\lambda$ and $p,q\in
  \bqd(\cF)$. Then $p\leq^* q$ if and only if there is $\vare<\lambda$ such
  that    
\[\big(\forall(\alpha,Z,F)\in q\big)\big(\forall A\in
F^+\big)\big(\alpha>\vare\ \Rightarrow\ \big(\exists (\alpha',Z',F')\in
p\big)\big(A\cap Z'\in (F')^+\big)\big).\]   
\item Let $p,q\in\bqd$. Then the following are equivalent:
\begin{enumerate}
\item[(a)]  $p\leq^* q$,
\item[(b)]  there is $\vare<\lambda$ such that   
\[\big(\forall(\alpha,Z,d)\in q\big)\big(\forall A\in
d\big)\big(\alpha>\vare\ \Rightarrow\ \big(\exists (\alpha',Z',d')\in p\big) 
\big(A\cap Z'\in d'\big)\big),\]    
\item[(c)] there is $\vare<\lambda$ such that if $(\alpha,Z,d)\in q$,
  $\vare\leq\alpha$, and $X=\big\{(\xi,Z',d')\in p: Z'\cap Z\neq\emptyset
  \big\}$, 
then $X\neq\emptyset$ and there is an ultrafilter $e$ on
$X$ such that  
\[d=\big\{A\cap Z:A\in\bigoplus\limits^e\{d':(\exists\xi,Z')( (\xi,Z',d')\in
X)\}\big\}.\]  
\end{enumerate}
\end{enumerate}
\end{proposition}

The quasi-orders $(\bqd,\leq^*)$ and $(\bqz,\leq^*)$ are
$({<}\lambda^+)$--complete (cf. \cite[Prop. 2.3(3)]{Sh:830}). Moreover, by
essentially the same argument we may show the following observation.
\begin{proposition}
  \label{whenclosed}
Assume that $\cF$ is a system of local filters on $\lambda$ such that 
\begin{enumerate}
\item[$(\oplus)^{\rm sum}_\cF$] if $\kappa<\lambda$ is an infinite cardinal
  and a sequence $\langle (\alpha_\xi,Z_\xi,F_\xi):\xi<\kappa\rangle
  \subseteq \cF$ satisfies 
\[(\forall\xi<\zeta<\kappa)(Z_\xi\subseteq\alpha_\zeta),\]
then for some uniform filter $F$ on $\kappa$ we have
$\big(\alpha_0,\bigcup\limits_{\xi<\kappa} Z_\xi,
\bigoplus\limits_{\xi<\kappa}^F F_\xi\big)\in \cF$.
\end{enumerate}
Then both $\bqd(\cF)$ and $\bqz(\cF)$ are $({<}\lambda^+)$--complete (with 
respect to $\leq^*$). 
\end{proposition}

It is worth noticing that in general $\bqd(\cF)$ and/or $\bqz(\cF)$ do not
have to be even $\sigma$--complete. For instance, consider the full system
of co-bounded filters $\cF_0$; it consists of all triples $(\alpha,Z,F)$
such that $\alpha\in Z\subseteq\lambda\setminus\alpha$, $|Z|<\lambda$,
$\sup(Z)\notin Z$ and $F$ is the filter of all co-bounded subsets of
$Z$. Let $C$ consist of all ordinals $\alpha<\lambda$ divisible by 
$\omega\cdot\omega$, and for $\alpha\in C$ and $m<\omega$ let
$Z^\alpha_m=[\alpha+m\cdot\omega,\alpha+m\cdot\omega+\omega)$ and 
$F^\alpha_m$ be the filter of co-bounded subsets of $Z^\alpha_m$. For $n<
\omega$ put 
\[p_n=\big\{(\alpha+m\cdot\omega,Z^\alpha_m,F^\alpha_m):\alpha\in C\ \&\
0<m<\omega\ \&\ 2^n|m\big\}.\]
Clearly $p_n\in\bqz(\cF_0)$ and $p_n\leq^* p_{n+1}$ for all $n<\omega$. One
may easily verify that the sequence $\langle p_n:n<\omega\rangle$ has no
$\leq^*$--upper bound in $\bqd(\cF_0)$.  

There is a natural procedure which for a given system $\cF$ of local filters
on $\lambda$ generates a system $\cF^*\supseteq \cF$ satisfying the
condition $(\oplus)^{\rm sum}_\cF$ of \ref{whenclosed} (so then
$\bqd(\cF^*)$ and $\bqz(\cF^*)$ are suitably complete).

\begin{definition}
 \label{sumclosure}
Assume that
\begin{enumerate}
\item[(a)] $\cF$ is a system of local filters on $\lambda$,
\item[(b)] $\bar{E}=\langle E_\kappa:\kappa\mbox{ is a cardinal }
\ \&\ \aleph_0\leq \kappa<\lambda\rangle$, where each $E_\kappa$ is a 
  uniform filter on $\kappa$.  
\end{enumerate}
We define:
\begin{enumerate}
\item An {\em $(\bar{E},\cF)$--block\/} is a pair $(T,\bar{D})$ such that 
  \begin{itemize}
  \item $T\subseteq{}^{<\omega}\lambda$ is a well--founded tree,
  \item if $\eta\in T\setminus \max(T)$, then $\{\xi<\lambda:\eta\conc
    \langle\xi\rangle\in T\}=\kappa$ for some infinite cardinal
    $\kappa<\lambda$,  
  \item $\bar{D}=\langle(\alpha_\eta,Z_\eta,F_\eta):\eta\in\max(T) \rangle
    \subseteq\cF$, 
  \item if $\eta,\nu\in\max(T)$ and $\eta<_{\rm lex}\nu$, then
    $Z_\eta\subseteq \alpha_\nu$ (where $<_{\rm lex}$ is the lexicographic
    order of $T$).
  \end{itemize}
\item By induction on the rank of the tree $T$, for an
  $(\bar{E},\cF)$--block $(T,\bar{D})$ we define a filter $\bar{D}(T)$ on
  $\bigcup\{Z_\eta:\eta\in \max(T)\}$ (where $\bar{D}=\langle( \alpha_\eta,
  Z_\eta,F_\eta):\eta\in\max(T) \rangle$).
  \begin{itemize}
  \item If $\rk(T)=0$, i.e., $T=\{\langle\rangle\}$ then $\bar{D}(T)=
    F_{\langle\rangle}$. 
  \item Suppose $\rk(T)>0$. Let $\kappa=\{\xi<\lambda:\langle\xi \rangle\in
    T\}$ (so $\aleph_0\leq\kappa<\lambda$ is a cardinal). For $\xi<\kappa$
    we put   
\[T^\xi=\{\nu\in {}^{<\omega}\lambda:\langle\xi\rangle\conc \nu\in T\}\
\mbox{ and }\ \bar{D}^\xi=\langle(\alpha_\eta,Z_\eta,F_\eta):\eta\in
\max(T)\ \&\ \eta(0)=\xi\rangle.\] 
Plainly, each $(T^\xi,\bar{D}^\xi)$ is an $(\bar{E},\cF)$--block (and
$\rk(T^\xi)<\rk(T)$). We define 
\[\bar{D}(T)=\bigoplus\limits_{\xi<\kappa}^{E_\kappa}\bar{D}^\xi(T^\xi).\]
  \end{itemize}
\item {\em The $\bar{E}$--closure of $\cF$\/} is the family of all triples
  $(\alpha,Z,D)$ such that $\alpha<\lambda$ and for some
  $(\bar{E},\cF)$--block $(T,\bar{D})$ we have 
\[Z=\bigcup\{Z_\eta:\eta\in\max(T)\}\quad\mbox{ and }\quad D=\bar{D}(T) 
\quad\mbox{ and }\quad \alpha=\min(Z)\]
(where $\bar{D}=\langle( \alpha_\eta, Z_\eta,F_\eta):\eta\in\max(T)
\rangle$). 
\end{enumerate}
\end{definition}

\begin{proposition}
  \label{closureclosed}
Assume that
\begin{enumerate}
\item[(a)] $\cF$ is a system of local filters on $\lambda$,
\item[(b)] $\bar{E}=\langle E_\kappa:\aleph_0\leq \kappa<\lambda\ \&\
  \kappa\mbox{ is a cardinal }\rangle$, where each $E_\kappa$ is a
  uniform filter on $\kappa$. 
\end{enumerate}
Then the $\bar{E}$--closure of $\cF$ is a system of local filters extending
$\cF$ and satisfying the condition $(\oplus)^{\rm sum}_\cF$ of
\ref{whenclosed}.  
\end{proposition}

Suppose that a system $\cF'$ of local filters on $\lambda$ includes all
triples $(\alpha,\{\alpha\},d)$, where $\alpha<\lambda$ and $d$ is the
principal ultrafilter on $\{\alpha\}$. For a set $A\subseteq\lambda$ let
$p_A=\{(\alpha,\{\alpha\},d)\in\cF': \alpha\in A\}\in \bqz(\cF')$. Note that
$p_A\cap p_B=p_{A\cap B}$, so easily if $D$ is a filter on $\lambda$
extending the co-bounded filter, then $H^D\stackrel{\rm def}{=} \{p_A:A\in
D\}$ is a $\leq^*$--directed family and $\fil(H^D)=D$. If $D$ is a normal
filter on $\lambda$, then $H^D$ will be also $({<}\lambda^+)$--directed (with
respect to $\leq^*$). Consequently, if $\lambda$ is a measurable cardinal,
then we may find a system $\cF$ of local filters on $\lambda$ and a
$({<}\lambda^+)$--directed family $H\subseteq \bqz(\cF)$ such that
$\fil(H)$ is an ultrafilter including all club subsets of $\lambda$ (so
$\fil(H)$ is not weakly reasonable). However, to have a quite directed
family $H$ such that $\fil(H)$ is a non-reasonable ultrafilter we do
need a measurable cardinal.

\begin{theorem}
  \label{nonres}
Suppose that $\cF$ is a system of local filters on $\lambda$,
$\kappa\leq\lambda$ and $H\subseteq\bqd(\cF)$ is a $({<}\kappa)$--directed
family such that $\fil(H)$ is an ultrafilter. If $\fil(H)$ is not
weakly reasonable, then for some club $C^*$ of $\lambda$ the quotient
ultrafilter $\fil(H)/C^*$ is $({<}\kappa)$--complete and it contains all
clubs of $\lambda$.
\end{theorem}

\begin{proof}
  Assume that the family $H\subseteq\bqd(\cF)$ is $({<}\kappa)$--directed
  and $\fil(H)$ is an ultrafilter which is not weakly reasonable. Let
  $\bar{\delta}=\langle\delta_\xi:\xi<\lambda\rangle$ be an increasing
  continuous sequence of ordinals below $\lambda$ such that $\delta_0=0$ and
  for every club $C\subseteq\lambda$ we have that $\bigcup\big\{[\delta_\xi,
  \delta_{\xi+1}):\xi\in C\big\}\in\fil(H)$. Now, for a club $C$ of
  $\lambda$ and $p\in H$ put  
\[S(p,C)=\big\{\xi\in C:\big(\exists (\alpha,Z,F)\in p\big)\big(
[\delta_\xi,\delta_{\xi+1})\cap Z\in F^+\big)\big\}.\]

\begin{claim}
\label{cl3}
For every club $C\subseteq\lambda$ and $p\in H$, the set $S(p,C)$ is
stationary.
\end{claim}

\begin{proof}[Proof of the Claim]
Assume towards contradiction that $S(p,C)$ is non-stationary. So we may
choose a club $C'\subseteq C$ of $\lambda$ such that 
\begin{enumerate}
\item[$(*)_1$] $\big(\forall\xi\in C'\big)\big(\forall (\alpha,Z,F)\in
  p\big)\big( Z\setminus [\delta_\xi,\delta_{\xi+1})\in F\big)$.
\end{enumerate}
Pick a club $C''\subseteq C'$ such that 
\begin{enumerate}
\item[$(*)_2$] $\big(\forall (\alpha,Z,F)\in p\big)\big(\forall\xi\in
  C''\big)\big(\alpha<\delta_\xi\ \Rightarrow\ Z\subseteq \delta_\xi\big)$.   
\end{enumerate}
By the choice of $\bar{\delta}$ we know that $\bigcup\{[\delta_\xi,
\delta_{\xi+1}):\xi\in C''\}\in\fil(H)$, so necessarily $\bigcup\{[\delta_\xi,
\delta_{\xi+1}):\xi\in C''\}\in\big(\fil(p)\big)^+$. Thus we may pick
$(\alpha,Z,F)\in p$ and $\xi\in C''$ such that $Z\cap
[\delta_\xi,\delta_{\xi+1})\in F^+$ (remember $(*)_2$), contradicting
$(*)_1$. 
\end{proof}

\begin{claim}
  \label{cl4}
  \begin{enumerate}
  \item If $p\leq^* q$, $p,q\in H$ and $C'\subseteq C$ are clubs of
    $\lambda$, then $|S(q,C')\setminus S(p,C)|<\lambda$.
  \item If $A\subseteq\lambda$, then there are $p\in H$ and a club
    $C\subseteq \lambda$ such that either $S(p,C)\subseteq A$ or
    $S(p,C)\subseteq \lambda\setminus A$. 
  \end{enumerate}
\end{claim}

\begin{proof}[Proof of the Claim]
  (1)\quad Pick $\gamma<\lambda$ so that 
\[\big(\forall (\alpha,Z,F)\in q\big)\big(\forall A\in F^+\big)\big(
\alpha>\gamma\ \Rightarrow\ (\exists (\alpha',Z',F')\in p)(A\cap Z'\in
(F')^+)\big)\] 
(remember \ref{1.3C}) and let $\gamma^*<\lambda$ be such that
$\gamma<\gamma^*$ and $\big(\forall (\alpha,Z,F)\in q\big)\big( \alpha\leq
\gamma\ \Rightarrow\ Z\subseteq\gamma^*\big)$. Suppose that $\xi\in
S(q,C')\setminus \gamma^*$. Then $\xi\in C'\subseteq C$ and there is
$(\alpha,Z,F)\in q$ such that $[\delta_\xi,\delta_{\xi+1})\cap Z\in
F^+$. Since $\delta_\xi\geq\xi\geq\gamma^*$, we also have $\alpha>\gamma$
and hence there is $(\alpha',Z',F')\in p$ such that $[\delta_\xi,
\delta_{\xi+1})\cap Z\cap Z'\in (F')^+$. Hence we may conclude that $\xi\in
S(p,C)$. 
\medskip

\noindent (2)\quad Assume $A\subseteq\lambda$. Let $A^*=\bigcup\big\{
[\delta_\xi,\delta_{\xi+1}):\xi\in A\big\}$. Since $\fil(H)$ is an
ultrafilter, then either $A^*$ or $\lambda\setminus A^*$ belongs to
it. Suppose $A^*\in\fil(p)$ for some $p\in H$. Pick a club
$C\subseteq\lambda$ such that  
\begin{enumerate}
\item[$(\odot)$] if $(\alpha,Z,F)\in p$ and $(\sup(Z)+1)\cap
  C\neq\emptyset$, then $A^*\cap Z\in F$. 
\end{enumerate}
Suppose $\xi\in S(p,C)$, so $\xi\in C$ and for some $(\alpha,Z,F)\in p$ we
have $[\delta_\xi,\delta_{\xi+1})\cap Z\in F^+$. It follows from $(\odot)$
that $A^*\cap Z\in F$ and therefore $\xi\in A$. Thus $S(p,C)\subseteq A$. 

If $\lambda\setminus A^*\in\fil(H)$, then we proceed in an analogous
manner. 
\end{proof}

Let 
\[D=\big\{A\subseteq\lambda:|S(p,C)\setminus A|<\lambda\mbox{ for some $p\in
  H$ and a club $C\subseteq\lambda$ }\big\}.\]
It follows from \ref{cl3} that all members of $D$ are stationary and since
$H$ is directed we may use \ref{cl4}(1) to argue that $D$ is a filter on
$\lambda$. By \ref{cl4}(2) we see that $D$ is an ultrafilter on $\lambda$
(so it also contains all clubs as its members are stationary). Since $H$ is
$({<}\kappa)$--directed and the intersection of ${<}\kappa$ many clubs is a
club, we may conclude from \ref{cl4}(1) that $D$ is a
$({<}\kappa)$--complete ultrafilter.

Let $C^*=\{\delta_\xi:\xi<\lambda\}$ (so it is a club of $\lambda$). To
complete the proof of the theorem we are going to show that $D=\fil(H)/C^*$.
Since we already know that $D$ is an ultrafilter, it is enough to show that
$S(p,C)\in \fil(H)/C^*$ for every $p\in H$ and a club $C\subseteq
\lambda$. So let $C\subseteq\lambda$ be a club, $p\in H$ and $S^*=
\bigcup\big\{[\delta_\xi,\delta_{\xi+1}): \xi\in S(p,C)\big\}$. If $S^*\in
\fil(H)$, then we are done, so assume that $S^*\notin\fil(H)$. Since
$\fil(H)$ is an ultrafilter and $H$ is directed, we may find $q\in H$ such
that $p\leq^* q$ and $\lambda\setminus S^*\in\fil(q)$. Let $\gamma<\lambda$
be such that
\[\big(\forall (\alpha,Z,F)\in q\big)\big(\gamma\leq\sup(Z)\ \Rightarrow\
Z\setminus S^*\in F\big).\]
Since $|S(q,C)\setminus S(p,C)|<\lambda$, we may pick $\xi\in S(q,C)\cap
S(p,C)$ such that $\xi>\gamma$. Then $[\delta_\xi,\delta_{\xi+1})\subseteq
S^*$ but also there is $(\alpha,Z,F)\in q$ such that
$[\delta_\xi,\delta_{\xi+1})\cap Z\in F^+$, and thus also $S^*\cap Z\in
F^+$. However, $\sup(Z)\geq\delta_\xi>\gamma$, so $Z\setminus S^*\in F$ by
the choice of $\gamma$, a contradiction showing that $S^*\in\fil(H)$ as
required. 
\end{proof}

\section{Systems of local ultrafilters}
In this section we are interested in the full system $\cF^{\rm ult}$ of
local ultrafilters on $\lambda$ and $\bqd,\bqz$. The first question that one
may ask is whether weakly reasonable ultrafilters on $\lambda$ generated by
some $H\subseteq\bqz(\cF)$ can be obtained by the use of $\bqz$. It occurs
that it does matter which system of local filters we are using.

\begin{definition}
  \label{UNultr}
A filter $F$ on a set $Z$ is called {\em an unultra filter, } if for every
$A\in F^+$ there is $B\subseteq A$ such that both $B\in F^+$ and $A\setminus 
B\in F^+$. The full system of local unultra filters on $\lambda$ will be
denoted by $\cF^{\rm unu}$. (Thus $\cF^{\rm unu}$ consists of all triples
$(\alpha, Z, F)$ such that $\emptyset\neq Z\subseteq \lambda$,
$|Z|<\lambda$, $\alpha=\min(Z)$ and $F$ is an unultra filter on $Z$.) 
\end{definition}

\begin{observation}
\label{unubasic}
\begin{enumerate}
\item If $F$ is an unultra filter on $Z$, $A\in F^+$, then $F+A
\stackrel{\rm def}{=}\{B\subseteq Z:B\cup (Z\setminus A)\in F\}$ is an
unultra filter.
\item Suppose that $\xi$ is a limit ordinal, $\{Z_\zeta:\zeta<\xi\}$ is a
  family of pairwise disjoint non-empty sets, $F_\zeta$ is a filter on
  $Z_\zeta$ (for $\zeta<\xi$). Then $\bigoplus\limits_{\zeta<\xi}F_\zeta$ is
  an unultra filter on $\bigcup\limits_{\zeta<\xi} Z_\xi$. (Remember the
  convention declared in the last sentence of Definition \ref{eproduct}.) 
\end{enumerate}
\end{observation}

\begin{theorem}
  \label{nonultra}
Assume $\lambda^{<\lambda}=\lambda$ and $2^\lambda=\lambda^+$. There exists
a $\leq^*$--increasing sequence $\langle p_\xi:\xi<\lambda^+\rangle\subseteq
\bqz(\cF^{\rm unu})$ such that  
\begin{enumerate}
\item[(a)]  $\fil\big(\{p_\xi:\xi<\lambda^+\}\big)$ is a weakly reasonable
  ultrafilter on $\lambda$, but 
\item[(b)] there is no $p\in\bqz$ with $\fil(p)\subseteq \fil
  \big(\{p_\xi: \xi<\lambda^+\}\big)$. 
\end{enumerate}
\end{theorem}

\begin{proof}
Fix enumerations 
\begin{itemize}
\item $\langle Y_\zeta:\zeta<\lambda^+\ \&\ \zeta$ is limit $\rangle$ of all
  subsets of $\lambda$, and
\item $\langle r_\zeta:\zeta<\lambda^+\ \&\ \zeta$ is limit $\rangle$ of
  $\bqz$, and 
\item $\langle\bar{\delta}^\zeta:\zeta<\lambda^+\ \&\ \zeta$ is limit
  $\rangle$ of all increasing continuous sequences of ordinals below 
$\lambda$, $\bar{\delta}^\zeta=\langle\delta^\zeta_\alpha:\alpha<\lambda
\rangle$.    
\end{itemize}
By induction on $\xi<\lambda^+$ we choose $p_\xi\in\bqz(\cF^{\rm unu})$ so
that the following conditions are satisfied for every limit  ordinal
$\zeta<\lambda^+$.  
\begin{enumerate}
\item[(o)]   For $n<\omega$, the element $p_n\in\bqz(\cF^{\rm unu})$ is 
\[\big\{(\alpha,Z_\alpha,F_\alpha):\alpha<\lambda\mbox{ is
    limit, }Z_\alpha=[\alpha,\alpha+\omega)\mbox{ and }F_\alpha\mbox{ is
    the filter of co-finite subsets of }Z_\alpha\big\}.\]
\item[(i)]  If $\cf(\zeta)<\lambda$, then for some increasing and cofinal in
$\zeta$ sequence $\langle\zeta_i:i<\cf(\zeta)\rangle$, for every
$(\alpha,Z,F)\in p_\zeta$, there is a sequence $\langle (\alpha_i,Z_i,F_i):
i<\cf(\zeta)\rangle$ such that 
\begin{itemize}
\item $(\alpha_i,Z_i,F_i)\in p_{\zeta_i}$,
\item $Z_i\subseteq \alpha_j$ for $i<j<\cf(\zeta)$,
\item $Z=\bigcup\limits_{i<\cf(\zeta)} Z_i$ and
  $F=\bigoplus\{F_i:i<\cf(\zeta)\}$. 
\end{itemize}
\item[(ii)]   If $\cf(\zeta)=\lambda$, then for some increasing and cofinal in
$\zeta$ sequence $\langle\zeta_i:i<\lambda\rangle$, if $(\alpha,Z,F)\in
p_\zeta$ and $\otp\big(\big\{\alpha'<\alpha:\big(\exists Z',F'\big)\big((\alpha',
Z', F')\in p_\zeta\big)\big\}\big)=j$, then  $(\alpha,Z,F)\in p_{\zeta_j}$
and 
\[\big(\forall i<j\big)\big(\forall A\in F^+\big)\big(\exists (\beta,W,D)\in
p_{\zeta_i}\big)\big(A\cap W\in D^+\big).\]
\item[(iii)] If $|\{(\alpha,Z,F)\in p_\zeta:Y_\zeta\cap Z\in F^+\}|=
\lambda$, then 
\[p_{\zeta+1}=\big\{\big(\alpha,Z,F+[Y_\zeta\cap Z]\big):
(\alpha,Z,F)\in p_\zeta\ \&\ Y_\zeta\cap Z\in F^+ \big\},\]
and otherwise $p_{\zeta+1}=\big\{(\alpha,Z,F)\in p_\zeta: Z\setminus
Y_\zeta\in F\big\}$. 
\item[(iv)] $p_{\zeta+2}\subseteq p_{\zeta+1}$ and for some club $C$ of
  $\lambda$, for every $\beta\in C$ we have 
\begin{itemize}
\item $Z\subseteq \delta^\zeta_\beta$ whenever $(\alpha,Z,F)\in
  p_{\zeta+2}$, $\alpha<\delta^\zeta_\beta$, and 
\item $\delta^\zeta_{\beta+1}<\min\big(\alpha\geq \delta^\zeta_\beta:
  \big(\exists Z,F\big)\big((\alpha,Z,F)\in p_{\zeta+2}\big)\big)$.   
\end{itemize}
\item[(v)] $p_{\zeta+3}=\big\{(\alpha,Z,F+A_\alpha):(\alpha,Z,F)\in
  p_{\zeta+2}\big\}$ where for every $(\alpha,Z,F)\in p_{\zeta+2}$ the set
 $A_\alpha\in F^+$ is such that $\big(\forall (\beta,Y,d)\in r_\zeta\big)
 \big(A_\alpha\cap Y\notin d\big)$.
\item[(vi)] $p_{\zeta+3+n}=p_{\zeta+3}$ for all $n<\omega$.
\item[(vii)] For every $r\in\bqz$ and $\xi<\lambda^+$, if $(\alpha,Z,F)\in 
  p_\xi$ and $A\in F^+$, then there is $A'\subseteq A$ such that   
\[A'\in F^+\quad \mbox{ and }\quad \big(\forall (\beta,Y,d)\in r\big)
 \big(A'\cap Y\notin d\big).\] 
\end{enumerate}
Conditions (o)--(vi) fully describe how the construction is carried out and
\ref{1.3C}+\ref{unubasic} imply that $\langle p_\xi:\xi<\lambda^+\rangle
\subseteq\bqz(\cF^{\rm unu})$ is $\leq^*$--increasing. However, we have to
argue that the demand in (vii) is satisfied, as it is crucial for the
possibility of satisfying the demand in (v). Let $r\in\bqz$. By induction on
$\xi<\lambda$ we show that for every $(\alpha,Z,F)\in p_\xi$ we have 
\begin{enumerate}
\item[$(\boxdot)_{(\alpha,Z,F)}$]  if  $A\in F^+$, then there
  is $A'\subseteq A$ such that  $A'\in F^+$ and $\big(\forall (\beta,Y,d)\in
  r\big)\big(A'\cap Y\notin d\big)$.
\end{enumerate}
(For a set $A'$ as above we will say that {\em it works for  $F$ and $r$.})
\medskip

\noindent {\sc Step} $\xi<\omega$.\\
Note that for each limit ordinal $\alpha<\lambda$ there is at most one
$(\beta,Y,d)\in r$ such that $Y\cap [\alpha,\alpha+\omega)$ is
infinite. Assume $A\subseteq [\alpha,\alpha+\omega)$ is
infinite. Considering any two disjoint infinite sets $A',A''\subseteq A$ we
easily see that one of them must work for the filter of co-finite subsets of
$[\alpha,\alpha+\omega)$ and $r$.  
\medskip

\noindent {\sc Step} $\xi=\zeta+n+1$, $\zeta<\lambda^+$ is limit,
$n<\omega$.\\ 
If $(\alpha,Z,F)\in p_{\zeta+n}$, $A\in F^+$, $A\subseteq A^*$ and
$A'\subseteq A$ works for $F$ and $r$, then also $A'$ works for $F+A^*$
and $r$.  
\medskip

\noindent {\sc Step} $\xi=\zeta<\lambda^+$ is limit.\\
Suppose that $(\alpha,Z,F)\in p_\zeta$. If $\cf(\zeta)=\lambda$, then
$(\alpha,Z,F)\in p_{\xi'}$ for some $\xi'<\zeta$ (see (ii)) so the inductive 
hypothesis applies directly. So assume that $\cf(\zeta)<\lambda$. Then 
$Z=\bigcup\limits_{i<\cf(\zeta)} Z_i$ and $F=\bigoplus\{F_i:i<\cf(\zeta)\}$
for some sequence $\langle (\alpha_i,Z_i,F_i):i<\cf(\zeta)\rangle$ such that 
\begin{itemize}
\item $(\boxdot)_{(\alpha_i,Z_i,F_i)}$ holds for each $i<\cf(\zeta)$, and 
\item $Z_i\subseteq \alpha_j$ for $i<j<\cf(\zeta)$.
\end{itemize}
Let $\alpha^*=\sup(Z)$  and let $A\in F^+$. Now we consider three cases. 
\smallskip

\noindent{\sc Case A:}\quad For some $\alpha'<\alpha^*$  we have
$\big(\forall (\beta,Y,d)\in r\big)\big(Y\cap
[\alpha',\alpha^*)=\emptyset\big)$.\\ 
Plainly, the set $A'=A\setminus \alpha'$ works for $F$ and $r$.
\smallskip

\noindent{\sc Case B:}\quad For some $(\beta,Y,d)\in r$ we have
$\beta<\alpha^*\leq \sup(Y)$. \\
For each $i<\cf(\zeta)$ such that $A\cap Z_i\in (F_i)^+$ choose disjoint
sets $A_i^0,A_i^1\in (F_i)^+$ included in $A\cap Z_i$ (remember each $F_i$ is
an unultra filter) and let 
\[A^\ell=\bigcup\{A_i^\ell:i<\cf(\zeta)\ \&\ A\cap Z_i\in
(F_i)^+\}\setminus\beta\subseteq A\] 
(for $\ell<2$). Both $A^0\in F^+$ and $A^1\in F^+$, and one of these two
sets works for $F$ and  $r$. 
\smallskip

\noindent{\sc Case C:}\quad For each $\alpha'<\alpha^*$ there is
$(\beta,Y,d)\in r$ such that $\alpha'<\beta<\sup(Y)<\alpha^*$. \\
Let $A\in F^+$. Then the set $I=\{i<\cf(\zeta):A\cap Z_i\in (F_i)^+\}$ is
unbounded in $\cf(\zeta)$ and using the assumptions of the current case we
may choose an increasing sequence $\langle i_j: j<\cf(\zeta)
\rangle\subseteq I$ such that for every $(\beta,Y,d)\in r$ there is at most
one $j<\cf(\zeta)$ such that  $Z_{i_j}\cap Y\neq\emptyset$. For 
each $j<\cf(\zeta)$ pick $A_{i_j}'\in (F_{i_j})^+$ included in $A\cap
Z_{i_j}$ which works for $F_{i_j}$ and $r$, and then put
$A'=\bigcup\limits_{j<\cf(\zeta)} A_{i_j}'$.  
\end{proof}

\begin{problem}
Is it provable in ZFC that for some system $\cF$ of local filters on
$\lambda$ there exists a $\leq^*$--directed family $H\subseteq\bqd(\cF)$ such that  
\begin{enumerate}
\item[(a)]  $\fil(H)$ is a weakly reasonable ultrafilter on $\lambda$, but 
\item[(b)]  there is no $\leq^*$--directed family $H'\subseteq \bqz(\cF^{\rm
    ult})$ such that $\fil(H)=\fil(H')$?     
\end{enumerate}
\end{problem}

The assumption that a generating system $H\subseteq\bqd(\cF)$ is directed
is an easy way to ensure that $\fil(H)$ is a filter on
$\lambda$. However, if we work with $H\subseteq\bqd$ we may consider
alternative ways of guaranteeing this. 

\begin{definition}
\label{bigA}
For $p\in\bqd$ let 
\[\Sigma(p)=\big\{(\alpha,Z,d)\in\cF^{\rm ult}: \big(\forall A\in d\big)
\big(\exists (\alpha',Z',d')\in p\big)\big(A\cap Z'\in d'\big)\big\}.\] 
\end{definition}

\begin{observation}
\label{bigB}
\begin{enumerate}
\item If $p,q\in\bqd$, then 
$p\leq^* q$\quad if and only if\quad $|q\setminus \Sigma(p)|<\lambda$. 
\item If $p\in\bqd$ and $(\alpha,Z,d)\in \Sigma(p)$, then for some
  $\{(\alpha_x,Z_x,d_x):x\in X\}\subseteq p$ and an ultrafilter $e$ on $X$
  we have 
\[d=\big\{A\subseteq Z:A\cap\bigcup\limits_{x\in X} Z_x\in\bigoplus\limits^e
  \{d_x:x\in X\}\big\}.\]
\end{enumerate}
\end{observation}

\begin{definition}
\label{bigC}
We say that a non-empty family $H\subseteq\bqd$ is 
\begin{enumerate}
\item[(a)] {\em big\/} if for each $\cD\subseteq \cF^{\rm ult}$ there is
  $q\in H$ such that either $q\subseteq \cD$ or $q\cap\cD=\emptyset$;   
\item[(b)] {\em linked\/} if for each $p_0,\ldots,p_n\in H$, $n<\omega$,
  we have 
\[\big\{\alpha:\big(\exists Z,d\big)\big((\alpha,Z,d)\in \Sigma(p_0)\cap
  \ldots\cap \Sigma(p_n)\big)\big\}|=\lambda.\] 
\end{enumerate}
\end{definition}

The property introduced in Definition \ref{bigC}(a) resembles the bigness of
creating pairs (see \cite[Sec. 2.2]{RoSh:470}), so the use of the term {\em
  big\/} seemed natural. The name {\em linked\/} is motivated by Observation
\ref{bigD}(1) below.    

\begin{observation}
\label{bigD}
\begin{enumerate}
\item If $H\subseteq\bqd$ is linked, then 
\begin{enumerate}
\item[(a)] for each $p_0,\ldots,p_n\in H$, $n<\omega$, there is $q\in\bqd$
  which is $\leq^*$--above all $p_0,\ldots,p_n$,
\item[(b)] $\fil(H)$ has finite intersection property.
\end{enumerate}
\item If $H\subseteq\bqd$ is linked and big, then $\fil(H)$ is an
  ultrafilter on $\lambda$.
\end{enumerate}
\end{observation}

For basic information on the ideal $\Mkkk$ of meager subsets of
${}^\lambda\lambda$ and its covering number we refer the reader e.g. to
Matet, Ros{\l}anowski and Shelah \cite[\S 4]{MRSh:799}. Let us recall the
following definition.

\begin{definition}
\label{d0.1}
\begin{enumerate}
\item The space ${}^\lambda\lambda$ is endowed with the topology obtained by
taking as basic open sets $\emptyset$ and $O_s$ for $s\in{}^{<\lambda}
\lambda$, where $O_s=\{f\in {}^\lambda\lambda: s\subseteq f\}$. 
\item The $({<}\lambda^+)$--complete ideal of subsets of ${}^\lambda\lambda$
  generated by nowhere dense subsets of ${}^\lambda\lambda$ is denoted by
  $\Mkkk$.
\item ${\rm cov}(\Mkkk)$ is the minimal size of a family $\cA\subseteq\Mkkk$
  such that $\bigcup\cA={}^\lambda\lambda$.
\end{enumerate}
\end{definition}

\begin{theorem}
\label{bigE}
Assume that $\lambda=\lambda^{<\lambda}\geq \aleph_1$ and ${\rm
  cov}(\Mkkk)=2^\lambda$. Then there exists a linked and big family
$H\subseteq\bqz$ such that $\fil(H)$ is a weakly reasonable ultrafilter.   
\end{theorem}

\begin{proof}
  The proof is very similar to that of \cite[Thm 2.14]{Sh:830}. Let $\chi$
  be a sufficiently large regular cardinal and let $N\prec\cHchi$ be such
  that $|N|=\lambda$ and ${}^{<\lambda}N\subseteq N$. Put $\cF_N^{\rm
    ult}=\cF^{\rm ult}\cap N$. We will inductively construct a linked and
  big family $H$ included in $\bqz(\cF_N^{\rm  ult})\subseteq\bqz(\cF^{\rm
    ult})$. The following two claims are the key points of the inductive
  process. Below, ``linked'' means ``linked as a subfamily of $\bqd$''
  (i.e., it is the notion introduced in Definition \ref{bigC}(b)).  

\begin{claim}
\label{cl1}
Assume that $H_0\subseteq\bqd(\cF_N^{\rm ult})$ is linked, $|H_0|<{\rm 
  cov}(\Mkkk)$, and $\cD\subseteq \cF^{\rm ult}$. Then there is
$q\in\bqz(\cF_N^{\rm ult})\subseteq \bqz$ such that $H_0\cup\{q\}$ is linked
and either $q\subseteq \cD$ or $q\cap\cD=\emptyset$.  
\end{claim}

\begin{proof}[Proof of the Claim]
We consider two cases.

\noindent {\sc Case A:}\quad For every $n<\omega$, $p_0,\ldots,p_n\in H_0$
and $\beta<\lambda$ there is $(\alpha,Z,d)\in \Sigma(p_0)\cap\ldots\cap 
\Sigma(p_n)\cap \cD\cap \cF_N^{\rm ult}$ such that $\beta<\alpha$. \\
Let $\cT_0$ be the family of all sequences $\eta$ such that 
\begin{enumerate}
\item[(i)]   $\lh(\eta)<\lambda$,
\item[(ii)]  if $\xi<\lh(\eta)$, then $\eta(\xi)\in \cD\cap \cF_N^{\rm
    ult}$,  
\item[(iii)] if $\xi<\xi'<\lh(\eta)$, $\eta(\xi)=(\alpha,Z,d)$, $\eta(\xi')
  =(\alpha',Z',d')$, then $Z\subseteq\alpha'$.
\end{enumerate}
It follows from the assumptions of the current case that $\cT_0$ is a
$\lambda$--branching tree (remember $|\cF_N^{\rm ult}|=\lambda$).  Moreover,
for each $p_0,\ldots,p_n\in H_0$ we have  
\[\big\{\rho\in\lim(\cT_0):\big(\exists\zeta<\lambda\big)\big(\forall\xi>
\zeta\big)\big(\rho(\xi)\notin \Sigma(p_0)\cap\ldots\cap \Sigma(p_n)\big)
\big\}\in\Mkkk.\]  
Hence (as $|H_0|<{\rm cov}(\Mkkk)$) we may pick $\rho\in\lim(\cT_0)$ such
that for every $p_0,\ldots,p_n\in H_0$, $n<\omega$, we have 
\[|\big\{\xi<\lambda:\rho(\xi)\in \Sigma(p_0)\cap\ldots\cap \Sigma(p_n)\big\}
|=\lambda.\] 
Let $q=\{\rho(\xi):\xi<\lambda\}$. Then $q\in \bqz(\cF_N^{\rm ult})\subseteq
\bqz$, $H_0\cup\{q\}$ is linked and $q\subseteq \cD$. 
\medskip

\noindent {\sc Case B:}\quad Not Case A.\\
Then for some $p^*_0,\ldots,p_m^*\in H_0$ and $\beta<\lambda$ we have  
\[\big(\forall (\alpha,Z,d)\in \Sigma(p_0^*)\cap\ldots\cap\Sigma(p^*_m) \cap
\cF_N^{\rm ult} \big)\big(\alpha>\beta\ \Rightarrow\ (\alpha,Z,d)\notin
\cD\big).\]  
It follows from the choice of $N$ that 

if $p_0,\ldots,p_n\in\bqd(\cF_N^{\rm ult})$ and $(\alpha,Z,d)\in\Sigma(p_0)
\cap \ldots\cap \Sigma(p_n)$,

then  there are $Z',d'$ such that $(\alpha,Z',d')\in\Sigma(p_0)\cap\ldots
\cap \Sigma(p_n) \cap \cF_N^{\rm ult}$. 

\noindent Consequently, we may repeat arguments of the previous case
replacing in clause (ii) $\cD\cap\cF_N^{\rm ult}$ by $\cF_N^{\rm ult}
\setminus \cD$. Then we obtain $q\in\bqz(\cF_N^{\rm ult})\subseteq\bqz$ such that
$H_0\cup\{q\}$ is linked and $q\cap\cD=\emptyset$. 
\end{proof}

\begin{claim}
\label{cl2}
Assume that $H_0\subseteq\bqd(\cF_N^{\rm ult})\subseteq\bqd$ is linked,
$|H_0|<{\rm cov}(\Mkkk)$ and a sequence $\langle\delta_\xi:\xi<\lambda
\rangle\subseteq \lambda$ is increasing continuously. Then there are
$p\in\bqz(\cF_N^{\rm ult})$ and a club $C^*$ of $\lambda$ such that   
\begin{enumerate}
\item[(a)] $H_0\cup\{p\}$ is linked, and 
\item[(b)] $\bigcup\big\{[\delta_{\xi+1},\delta_\zeta):\xi<\zeta$ are
  successive members of $ C^*\big\}\in \fil(p)$.  
\end{enumerate}
\end{claim}

\begin{proof}[Proof of the Claim] 
This is essentially \cite[Claim 2.14.4]{Sh:830}.
\end{proof}

Now we employ a bookkeeping device to construct inductively a sequence
$\langle q_\xi:\xi<2^\lambda\rangle\subseteq\bqz(\cF_N^{\rm ult})$ such that 
\begin{itemize}
\item for each $\zeta<2^\lambda$ the family $\{q_\xi:\xi<\zeta\}$ is linked, 
\item if $\cD\subseteq \cF_N^{\rm ult}$, then for some $\xi<2^\lambda$ we
  have $q_\xi\subseteq \cD$ or $q_\xi\cap \cD=\emptyset$, 
\item if $\langle\delta_\xi:\xi<\lambda\rangle\subseteq\lambda$ is
  increasing continuous, then for some $\vare<2^\lambda$ and a club $C^*$ of 
  $\lambda$ we have that 
\[\bigcup\big\{[\delta_{\xi+1},\delta_\zeta):\xi
  <\zeta\mbox{ are successive members of } C^*\big\}\in \fil(q_\vare).\]  
\end{itemize}
Since $|\cF_N^{\rm ult}|=\lambda$, so there are no problems with carrying
out the construction. It should be clear that at the end the family
$\{q_\xi:\xi<2^\lambda\}$ is linked, big and it generates a weakly
reasonable ultrafilter.
\end{proof}

Note that we may modify the construction in the proof of Theorem \ref{bigE}
so that the resulting $H$ is directed. Namely, by an argument similar to the
one in Claim \ref{cl1} we may show, that if $H_0\subseteq \bqd(\cF_N^{\rm
  ult})$ is linked, $|H_0|<{\rm cov}(\Mkkk)$ and $p_0,p_1\in H_0$, then
there is $q\in \bqz(\cF_N^{\rm ult})$ such that $q\subseteq
\Sigma(p_0)\cap\Sigma(p_1)$ and $H_0\cup\{q\}$ is linked. With this claim in
hands we may modify the inductive choice of $\langle q_\xi:\xi<
2^\lambda\rangle$ so that at the end $\{q_\xi:\xi<2^\lambda\}$ is
directed. However, we do not know how to guarantee the opposite, that the
family $\{q_\xi:\xi<2^\lambda\}$ is not directed or even better, that for
{\bf no} directed $H\subseteq \bqz$ do we have
$\fil(H)=\fil(\{q_\xi:\xi<2^\lambda\})$. Thus the following question remains
open. 

\begin{problem}
\label{prbig}
Does ``$H\subseteq\bqz$ is linked and big'' imply that ``$H$ is directed''? 
\end{problem}

\section{Systems of local pararegular filters}
In this section we are interested in filters associated with the full system
$\cF^{\pr}$ of local pararegular filters on $\lambda$ and we show their
relation to numbers of generators (in standard sense) of some filters on
$\lambda$.

\begin{definition}
\label{spr}
Suppose that $Z\subseteq\lambda$ is an infinite set,  $\alpha=\min(Z)$. A
{\em pararegular filter on $Z$} is a filter $F$ on $Z$ such that for some
system $\langle A_u:u\in [\kappa]^{<\omega}\rangle$ of sets from $F$ we
have:  
\begin{itemize}
\item $|\omega+\alpha|\leq\kappa<\lambda$, and if $u\subseteq v\in
  [\kappa]^{<\omega}$, then $A_v\subseteq A_u$,  
\item if $U\subseteq\kappa$ is infinite, then $\bigcap\{A_{\{\xi\}}:\xi \in
  U\}=\emptyset$, and  
\item $\cF=\big\{B\subseteq Z:\big(\exists u\in [\kappa]^{<\omega}\big)
  \big(A_u\subseteq B\big)\big\}$.
    \end{itemize}
If the cardinal $\kappa$ above satisfies $2^{|\omega+\alpha|}\leq\kappa
<\lambda$, then we say that the filter $F$ is {\em strongly pararegular.}
 
The full system of local pararegular filters on $\lambda$ will be
denoted by $\cF^{\pr}$ and the full system of local strongly pararegular
filters on $\lambda$ is denoted by $\cF^{\rm spr}$. (The latter forms a
system of local filters if and only if $\lambda$ is inaccessible and then
$\cF^{\rm spr}\subseteq\cF^{\pr}$.) 
\end{definition}

Let us recall the following strong $\lambda^+$--chain condition. 

\begin{definition}
[See Shelah {\cite[Def. 1.1]{Sh:288}} and {\cite[Def. 7]{Sh:546}}]  
\label{288cc}
Let $\bbQ$ be a forcing notion, and $\varepsilon<\lambda$ be a limit ordinal. 
\begin{enumerate}
\item We define a game $\Gcc(\bbQ)$ of two players, Player I and Player II. A
play lasts $\varepsilon$ steps, and at each stage $\alpha<\varepsilon$ of
the play sequences $\bar{p}^\alpha,\bar{q}^\alpha$ and a function
$\varphi^\alpha$ are chosen so that: 
\begin{itemize}
\item $\bar{p}^0=\langle \emptyset_\bbQ:i<\lambda^+\rangle$, $\varphi^0:
\lambda^+\longrightarrow\lambda^+:i\mapsto 0$;
\item If $\alpha>0$, then Player I picks $\bar{p}^\alpha,\varphi^\alpha$
such that   
\begin{enumerate}
\item[(i)] $\bar{p}^\alpha=\langle p^\alpha_i:i<\lambda^+\rangle\subseteq
\bbQ$ satisfies $(\forall \beta<\alpha)(\forall i<\lambda^+)(q^\beta_i\leq
p^\alpha_i)$, 
\item[(ii)] $\varphi^\alpha:\lambda^+\longrightarrow\lambda^+$ is regressive,
i.e., $(\forall i<\lambda^+)(\varphi^\alpha(i)<1+i)$; 
\end{enumerate}
\item Player II answers choosing a sequence $\bar{q}^\alpha=\langle
q^\alpha_i: i<\lambda^+\rangle\subseteq\bbQ$ such that $(\forall i<\lambda^+)
(p^\alpha_i\leq q^\alpha_i)$.
\end{itemize}
If at some stage of the game Player I does not have any legal move, then he
loses. If the game lasted $\varepsilon$ steps, Player I wins a play
$\langle\bar{p}^\alpha,\bar{q}^\alpha,\varphi^\alpha:\alpha<\varepsilon
\rangle$ if there is a club $C$ of $\lambda^+$ such that for each distinct
members $i,j$ of $C$ satisfying $\cf(i)=\cf(j)=\lambda$ and
$(\forall\alpha<\varepsilon)(\varphi^\alpha(i)=\varphi^\alpha(j))$, the set
$\{q^\alpha_i:\alpha<\varepsilon\}\cup \{q^\alpha_j: \alpha<\varepsilon\}$
has an upper bound in $\bbQ$. 

\item The forcing notion $\bbQ$ satisfies condition $(*)^\varepsilon_\lambda$
if Player I has a winning strategy in the game $\Gcc(\bbQ)$.
\end{enumerate}
\end{definition}

\begin{proposition}
[See Shelah {\cite[Iteration Lemma 1.3]{Sh:288}} and {\cite[Thm
      35]{Sh:546}}]  
\label{288it}
Let $\varepsilon<\lambda$ be a limit ordinal,
$\lambda=\lambda^{<\lambda}$. Suppose that $\bar{\bbQ}=\langle\bbP_\xi,\dbQ_\xi:
\xi<\gamma\rangle$ is a $({<}\lambda)$--support iteration such that for each
$\xi<\gamma$  
\[\forces_{\bbP_\xi}\mbox{`` }\dbQ_\xi\mbox{ satisfies
}(*)^\varepsilon_\lambda\mbox{ ''.}\]
Then $\bbP_\gamma$ satisfies $(*)^\varepsilon_\lambda$.
\end{proposition}

\begin{definition}
  \label{domfor}
Suppose that $D$ is a uniform filter on $\lambda$. We define a forcing
notion $\bbQ^{\pr}_D$ by:

\noindent {\bf a condition} is a tuple $p=\big(\zeta^p,\langle\alpha^p_\xi:
\xi\leq\zeta^p\rangle, \langle Z^p_\xi,F^p_\xi:\xi<\zeta^p\rangle,\cA^p
\big)$ such that 
\begin{enumerate}
\item[$(\alpha)$] $\cA^p\subseteq D$, $|\cA^p|<\lambda$, $\zeta^p<\lambda$, 
\item[$(\beta)$]   $\langle\alpha^p_\xi:\xi\leq\zeta^p\rangle$ is an
  increasing continuous sequence of ordinals below $\lambda$, 
\item[$(\gamma)$]  $Z^p_\xi=[\alpha^p_\xi,\alpha^p_{\xi+1})$ and $F^p_\xi$
  is a pararegular filter on $Z^p_\xi$;
\end{enumerate}

\noindent {\bf the order } $\leq_{\bbQ^{\pr}_D}=\leq$ is given by:\qquad
$p\leq q$ if and only if ($p,q\in \bbQ^{\pr}_D$ and)
\begin{enumerate}
\item[(i)]   $\cA^p\subseteq \cA^q$, $\zeta^p\leq \zeta^q$, 
\item[(ii)]  $\alpha^q_\xi=\alpha^p_\xi$ for $\xi\leq\zeta^p$, and
  $Z^p_\xi=Z^q_\xi$, $F^q_\xi=F^p_\xi$ for $\xi<\zeta^p$, 
\item[(iii)] if $A\in\cA^p$ and $\zeta^p\leq\xi<\zeta^q$, then $A\cap
  Z^q_\xi\in F^q_\xi$.
\end{enumerate}
\end{definition}

\begin{proposition}
 \label{easyfor}
Assume $\lambda^{<\lambda}=\lambda$ and let $D$ be a uniform filter on
$\lambda$. Then:
\begin{enumerate}
\item $\bbQ^{\pr}_D$ is a $({<}\lambda)$--complete forcing notion of size
  $2^\lambda$,
\item $\bbQ^{\pr}_D$ satisfies the condition $(*)^\vare_\lambda$ of
  \ref{288cc} for each limit ordinal $\vare<\lambda$,  
\item if $\name{r}$ is a $\bbQ^{\pr}_D$--name such that 
\[\forces_{\bbQ^{\pr}_D}\name{r}=\big\{(\alpha^p_\xi,Z^p_\xi,F^p_\xi):\xi<\zeta^p
\ \&\ p\in \Gamma_{\bbQ^{\pr}_D}\big\},\]
then $\forces_{\bbQ^{\pr}_D}$`` $\name{r}\in \bqz(\cF^{\pr})$ and
$D\subseteq\fil(\name{r})$ ''.  
\end{enumerate}
\end{proposition}

\begin{proof}
(1)\quad Note that if $\alpha<\beta<\lambda$, then there are
$\leq\sum\limits_{\kappa<\lambda}2^{\kappa\cdot |[\alpha,\beta)|}$ many
pararegular filters on $[\alpha,\beta)$. Hence easily
$|\bbQ^{\pr}_D|=2^\lambda$.  

If $\langle p_\alpha:\alpha<\gamma\rangle\subseteq \bbQ^{\pr}_D$ is
$\leq_{\bbQ^{\pr}_D}$--increasing, $\gamma<\lambda$, then letting
$\cA^q=\bigcup\limits_{\alpha<\gamma}\cA^{p_\alpha}$, 
\[\zeta^q=
\sup(\zeta^{p_\alpha}:\alpha<\gamma)\ \mbox{ and }\ \langle \alpha^q_\xi, 
Z^q_\xi,F^q_\xi:\xi<\zeta^q\rangle=\bigcup_{\alpha<\gamma}\langle
\alpha^q_\xi, Z^{p_\alpha}_\xi,F^{p_\alpha}_\xi:\xi<\zeta^{p_\alpha}\rangle\]   
we get a condition $q=(\zeta^q,\langle\alpha^q_\xi:\xi\leq\zeta^q\rangle,
\langle Z^q_\xi,F^q_\xi:\xi<\zeta^q\rangle,\cA^q)\in \bbQ^{\pr}_D$ stronger
than all $p_\alpha$ (for $\alpha<\gamma$). 
\medskip

\noindent (2)\quad Let $\cX$ consists of all sequences $\langle Z_\xi,F_\xi:
\xi<\zeta\rangle$ such that $\langle Z_\xi,F_\xi:\xi<\zeta\rangle=\langle
Z_\xi^p,F_\xi^p: \xi<\zeta^p\rangle$ for some $p\in \bbQ^{\pr}_D$. By what
we said earlier, $|\cX|=\lambda$, so we may fix an enumeration $\langle
\bar{\sigma}_\alpha:\alpha<\lambda\rangle$ of $\cX$. Now, let $\st$ be a
strategy of Player I in $\Gcc(\bbQ^{\pr}_D)$ which, at a stage
$\alpha<\vare$ of the play, instructs her to choose a legal inning
$\bar{p}^\alpha,\varphi^\alpha$ such that if $\lambda\leq i<\lambda^+$, then
$\langle Z_\xi^{p^\alpha_i}, F_\xi^{p^\alpha_i}:\xi<\zeta^{p^\alpha_i}
\rangle= \bar{\sigma}_{\varphi^\alpha(i)}$. (Note that there are legal
innings for Player I by the completeness of the forcing proved in (1)
above.) Plainly, if $\langle\bar{p}^\alpha,\bar{q}^\alpha, \varphi^\alpha:
\alpha<\varepsilon\rangle$ is a play  of $\Gcc(\bbQ^{\pr}_D)$ in which
Player I follows $\st$ and $\lambda\leq i<j<\lambda^+$ are such that
$\varphi^\alpha(i)=\varphi^\alpha(j)$ for all $\alpha<\vare$, then the
family $\{p^\alpha_i,p^\alpha_j:\alpha<\vare\}$ has an upper bound. Thus 
$\st$ is a winning strategy for Player I.
\medskip

\noindent (3)\quad Suppose $p\in \bbQ^{\pr}_D$ and let $\kappa= |\cA^p|
+|\omega+\alpha^p_{\zeta^p}|$. Fix a sequence $\langle A_\beta:\beta<
\kappa\rangle$ listing all members of $\cA^p\cup\{\lambda\}$ (with possible
repetitions) and let  $\langle u_\gamma:\gamma<\kappa\rangle$ be an
enumeration of $[\kappa]^{<\omega}$. By induction on $\gamma<\kappa$ choose
an increasing sequence $\langle \xi_\gamma:\gamma<\kappa\rangle\subseteq
[\alpha^p_{\zeta^p},\lambda)$ such that $\xi_\gamma\in
\bigcap\limits_{\beta\in u_\gamma} A_\beta$. (Remember, $D$ is a uniform
filter and $\kappa<\lambda$.) Let $\zeta^q=\zeta^p+1$, $\alpha^q_{\zeta^q} =
\sup(\xi_\gamma: \gamma<\kappa)+1$ and  for $u\in [\kappa]^{<\omega}$ let
$B_u=\{\xi_\gamma: u\subseteq u_\gamma\ \&\ \gamma<\lambda\}$. Then
$F^q_{\zeta^p}=\{B\subseteq [\alpha^q_{\zeta^p},\alpha^q_{\zeta^q}):
(\exists u\in [\kappa]^{<\omega})(B_u\subseteq B)\}$ is a pararegular filter
on $[\alpha^q_{\zeta^p},\alpha^q_{\zeta^q})$ and $A\cap
[\alpha^q_{\zeta^p},\alpha^q_{\zeta^q})\in F^q_{\zeta^p}$ for all
$A\in\cA^q$.  So now we may take a condition $q\in \bbQ^{\pr}_D$ stronger 
than $p$ and such that $Z^q_{\zeta^p}=[\zeta^p,\zeta^q)$,
$\cA^p=\cA^q$. Then $q\forces (\alpha^q_{\zeta^p},Z^q_{\zeta^p},
F^q_{\zeta^p}) \in \name{r}$. 

So we easily conclude that indeed $\forces_{\bbQ^{\pr}_D}$`` $\name{r}\in
\bqz(\cF^{\pr})$ and $D\subseteq\fil(\name{r})$ '' (remember the definition
of the order on $\bbQ^{\pr}_D$, specifically \ref{domfor}(iii)).  
\end{proof}

\begin{corollary}
  \label{iterdom}
Assume $\lambda^{<\lambda}=\lambda$, $2^\lambda=\lambda^+$,
$2^{\lambda^+}=\lambda^{++}$. Then there is a $({<}\lambda)$--complete   
$\lambda^+$--cc forcing notion $\bbP$ such that 
\[
\begin{array}{ll}
\forces_{\bbP}&\mbox{`` }2^\lambda=\lambda^{++}\mbox{ and if $D$ is a
  uniform filter on }\lambda\mbox{ generated by }\\
&\mbox{ less than $\lambda^{++}$ elements, then $D\subseteq
  \fil(r)$ for some $r\in\bqz(\cF^{\pr})$ ''.}
\end{array}\]
\end{corollary}

\begin{proof}
Using a standard bookkeeping argument build a ${<}\lambda$--support
iteration $\bar{\bbQ}=\langle\bbP_\xi,\dbQ_\xi:\xi<\lambda^{++}\rangle$ such
that 
\begin{itemize}
\item for each  $\xi<\lambda^{++}$ we have that $\forces_{\bbP_\xi}$``
  $\dbQ_\xi=\bbQ^{\pr}_{\name{D}}$ '' for some $\bbP_\xi$--name $\name{D}$
  for a uniform filter on $\lambda$,  
\item if $\langle \name{A}_\beta:\beta<\lambda^+\rangle$ is a sequence of
  $\bbP_{\lambda^{++}}$--names for subsets of $\lambda$, then for some
  $\xi<\lambda^{++}$ such that every $\name{A}_\beta$ is a $\bbP_\xi$--name
  we have 
\[\forces_{\bbP_\xi}\mbox{`` if }\langle
  \name{A}_\beta:\beta<\lambda^+\rangle\mbox{ generates a uniform filter $D$ on
  $\lambda$, then }\dbQ_\xi=\bbQ^{\pr}_{D}\mbox{ ''.}\]
\end{itemize}
Now look at the limit $\bbP_{\lambda^{++}}=\lim(\bar{\bbQ})$ (and remember
\ref{easyfor}, \ref{288it}). 
\end{proof}

\begin{proposition}
  \label{notcovered}
Assume $2^\lambda=\lambda^+$. Then there is a uniform ultrafilter $D$ on
$\lambda$ containing no $\fil(p)$ for $p\in \bqd(\cF^{\pr})$.
\end{proposition}

\begin{proof}
First note that if $F$ is a pararegular filter on $Z$, then for each
$\beta$ we have $Z\setminus\{\beta\}\in F$. Consequently, if $\cA \subseteq
[\lambda]^\lambda$ is a family with fip, $\{[\alpha,\lambda):\alpha<
\lambda\}  \subseteq\cA$, $|\cA|\leq\lambda$, and $p\in \bqd(\cF^{\pr})$,
then we may choose $A\subseteq\lambda$ such that 
\begin{itemize}
\item $\cA\cup\{A\}$ has fip,
\item for each $(\alpha,Z,F)\in p$ we have $|Z\cap A|\leq 1$ so also
  $Z\setminus A\in F$. 
\end{itemize}
Hence, by induction on $\xi<\lambda^+$, we may choose a sequence $\langle
A_\xi:\xi<\lambda^+\rangle$ of unbounded subsets of $\lambda$ such that 
\begin{itemize}
\item for $\xi<\lambda$, $A_\xi=[\xi,\lambda)$,
\item $\{A_\xi:\xi<\lambda^+\}$ has fip, 
\item for every $A\subseteq \lambda$ there is $\xi<\lambda^+$ such that
  either $A_\xi\subseteq A$ or $A_\xi\cap A=\emptyset$,
\item for every  $p\in \bqd(\cF^{\pr})$ there is $\xi<\lambda^+$ such that
  $\lambda\setminus A_\xi\in\fil(p)$. 
\end{itemize}
Then $D=\{A\subseteq\lambda:A_{\xi_0}\cap\ldots\cap A_{\xi_n}\subseteq A$
for some $\xi_0,\ldots,\xi_n<\lambda^+$, $n<\omega\}$ is an ultrafilter as
required.  
\end{proof}

\begin{proposition}
 \label{criterion}
Assume that 
\begin{enumerate}
\item[(a)] there exists a $\lambda$--Kurepa tree with $2^\lambda$
  $\lambda$--branches, 
\item[(b)] $D$ is a uniform filter on $\lambda$, 
\item[(c)] $p\in\bqz(\cF^{\pr})$ is such that $\fil(p)\subseteq D$, 
\item[(d)] if $\lambda$ is a limit cardinal, then it is strongly
  inaccessible and  $p\in\bqz(\cF^{\rm spr})$. 
\end{enumerate}
Then the filter $D$ cannot be generated by less than $2^\lambda$ sets, i.e.,
for every family $\cX\subseteq D$ of size less than $2^\lambda$ there is a
set $A\in D$ such that $|X\setminus A|=\lambda$ for all $X\in\cX$. 
\end{proposition}

\begin{proof}
Let $T$ be a $\lambda$--Kurepa tree with $2^\lambda$  $\lambda$--branches
(so each level in $T$ is of size $<\lambda$). For $\xi<\lambda$ let $T_\xi$
be the $\xi^{\rm th}$ level of $T$. Choose an increasing continuous
sequence $\langle\alpha_\xi:\xi<\lambda\rangle$ such that if
$(\alpha,Z,F)\in p$ and $\alpha_\xi\leq\alpha<\alpha_{\xi+1}$, then 
\begin{itemize}
\item $Z\subseteq\alpha_{\xi+1}$ and 
\item there is a  system $\langle A^\alpha_u:u\in [\kappa_\alpha]^{<\omega}
  \rangle$ of sets from $F$ witnessing that $F$ is pararegular (strongly
  pararegular if $\lambda$ is inaccessible) with $\kappa_\alpha$ satisfying
  $|T_\xi|\leq \kappa_\alpha$.  
\end{itemize}
For each $\xi<\lambda$ and $(\alpha,Z,F)\in p$ such that $\alpha_\xi \leq
\alpha<\alpha_{\xi+1}$, let us fix an injection $\pi^\alpha_\xi:T_\xi
\stackrel{1-1}{\longrightarrow} \kappa_\alpha$, and next for every 
$\lambda$--branch $\eta$ through $T$ let us choose a set 
$A_\eta\in D$ so that 
\begin{itemize}
\item if $\xi<\lambda$, $\nu\in T_\xi\cap\eta$, $(\alpha,Z,F)\in p$,
  $\alpha_\xi\leq\alpha<\alpha_{\xi+1}$, then $A_\eta\cap
  Z=A^\alpha_{\{\pi^\alpha_\xi(\nu)\}}$.  
\end{itemize}
For our conclusion, it is enough to show that if $B\in D$, then there are at
most finitely many $\lambda$--branches $\eta$ through $T$ such that
$|B\setminus A_\eta|<\lambda$. So suppose towards contradiction
$\eta_0,\eta_1,\eta_2,\ldots$  are distinct $\lambda$--branches through $T$,
$B\in D$ and $|B\setminus A_{\eta_n}|<\lambda$ for each $n<\omega$. The set
$\{(\alpha,Z,F)\in p:B\cap Z\in F^+\}$ is of cardinality $\lambda$, so we may 
find $\xi<\lambda$ and $\nu_n\in T_\xi$ (for $n<\omega$) such that 
\begin{itemize}
\item $\eta_n\cap T_\xi=\{\nu_n\}$ and $\nu_n\neq \nu_m$ for
  distinct $n,m$, and  
\item $B\cap Z^*\in (F^*)^+$ for some $(\alpha^*,Z^*,F^*)\in p$ satisfying 
  $\alpha_\xi\leq \alpha^*<\alpha_{\xi+1}$, and  
\item $B\setminus \alpha_\xi\subseteq A_{\eta_n}$ for all $n<\omega$.
\end{itemize}
Then $\emptyset\neq B\cap Z^*\subseteq
\bigcap\{A^\alpha_{\{\pi^\alpha_\xi(\nu_n)\}}:n<\omega\}$, a contradiction.
\end{proof}

\section{Forcing a very reasonable ultrafilter}
Our goal here is to show that the inaccessibility of $\lambda$ in the
assumptions of \cite[Prop. 1.6(1)]{Sh:830} is needed. This answers the
request of the referee of \cite{Sh:830} and fulfills the promise stated in
\cite[Rem. 1.7]{Sh:830}. Assuming that $\kappa$ is strongly inaccessible,
we will construct a CS iteration $\langle \bbP_\alpha,\name{\bbQ}_\alpha:
\alpha<\kappa\rangle$ of proper forcing notions such that 
\[\begin{array}{ll}
\forces_{\bbP_\kappa}&\mbox{`` there is a $(\leq\omega_1)$--directed family
  $H\subseteq \bqo$ such that }\\
&\ \fil(H)\mbox{ is a weakly reasonable ultrafilter on $\omega_1$ and yet
  }\\
&\mbox{ Odd has a winning strategy in  $\Game_{\fil(H)}$ ''.}
  \end{array}\] 
Let us recall the following definition.

\begin{definition}
[Shelah {\cite[Def. 1.4]{Sh:830}}] 
\label{1.5X}
Let $D$ be a uniform ultrafilter on $\lambda$. We define a game $\Game_D$
between two players, Odd and Even, as follows. A play of $\Game_D$ lasts
$\lambda$ steps and during a play an increasing continuous sequence
$\bar{\alpha}=\langle\alpha_i:i<\lambda \rangle\subseteq\lambda$ is
constructed. The terms of $\bar{\alpha}$ are chosen successively by the two
players so that Even chooses the $\alpha_i$ for even $i$ (including limit
stages $i$ where she has no free choice) and Odd chooses $\alpha_i$ for odd
$i$. Even wins the play if and only if $\bigcup\{[\alpha_{2i+1},
\alpha_{2i+2}):i<\lambda\}\in D$. 
\end{definition}

The following result was shown in \cite[Prop. 1.6]{Sh:830}:

\begin{proposition}
\label{1.7}
Assume $D$ is a uniform ultrafilter on $\lambda$.
\begin{enumerate}
\item If $\lambda$ is strongly inaccessible and Odd has a winning strategy
  in $\Game_D$, then $D$ is not weakly reasonable.
\item If $D$ is not weakly reasonable, then Odd has a winning strategy in
  the game $\Game_D$.
\end{enumerate}
\end{proposition}

Before we define our CS iteration $\langle \bbP_\alpha,\name{\bbQ}_\alpha:
\alpha<\kappa\rangle$ let us introduce two main ingredients used in the 
construction. 
\medskip

\noindent{\bf Sealing the branches:}\quad At each stage of the iteration we 
will first use forcing notions introduced in Shelah \cite[Ch. XVII, \S 
  2]{Sh:f}. 

For a tree $T\subseteq {}^{<\omega_1}\omega_1$, the set of all
$\omega_1$-branches through $T$ will be denoted by $\lim(T)$. Thus
$\lim(T)=\{\eta\in {}^{\omega_1}\omega_1:(\forall\alpha<\omega_1)(\eta
\rest\alpha\in T)\}$. 

\begin{lemma}
[Shelah {\cite[Ch. XVII, Fact 2.2]{Sh:f}}]
\label{CLpres}
Suppose that $T\subseteq {}^{<\omega_1}\omega_1$ is a tree of
height $\omega_1$. Let $\bbC$ be the Cohen forcing and $\name{\bbL}$ be
a $\bbC$--name for the Levy collapse of $2^{\aleph_2}$ to $\aleph_1$ (with 
countable conditions, so it is a $\sigma$--closed forcing notion). Then\quad  
$\forces_{\bbC*\name{\bbL}}\mbox{`` }\lim(T)=\big(\lim(T)\big)^\bV \mbox{
  ''}$.
\end{lemma}

\begin{definition}
[Shelah {\cite[Ch. XVII, Def. 2.3]{Sh:f}}]
\label{sealdef}
Suppose that $T\subseteq {}^{<\omega_1}\omega_1$ is a tree of
height $\omega_1$, $|T|=\aleph_1$, $|\lim(T)|\leq\aleph_1$. Let $\langle
B_i:i<\omega_1\rangle$ list all members of $\lim(T)$ (with possible
repetitions) and $\langle y_i:i<\omega_1\rangle$ list all elements of $T$ so
that [$y_j\vtl y_i\ \Rightarrow\ j<i]$. For $j<\omega_1$ we define  
\[B^*_j=\left\{\begin{array}{ll}
B_i     &\mbox{ if }j=2i,\\
\{y_i\} &\mbox{ if }j=2i+1,
	       \end{array}\right. \qquad\mbox{ and }\qquad 
B'_j=B^*_j\setminus \bigcup_{i<j}B^*_i.\]
Let $w=\{j<\omega_1:B_j'\neq\emptyset\}$ and for $j\in w$ let
$x_j=\min(B_j')$. Finally, we put $A=\{x_i:i\in w\}$. We define a forcing
notion $\bbP_T$ for {\em sealing the branches of $T$\/}:\\
{\bf a condition $p$ in }$\bbP_T$ is a finite function from $\dom(p)\subseteq
  A$ into $\omega$ such that if $\rho,\nu\in\dom(p)$ and $\rho\vtl\nu$, then
  $p(\eta)\neq p(\nu)$,\\
{\bf the order $\leq_{\bbP_T}$ of $\bbP_T$} is the inclusion, i.e., $p\leq
  q$ if and only if ($p,q\in\bbP_T$ and) $p\subseteq q$.
\end{definition}  

\begin{lemma}
[Shelah {\cite[Ch. XVII, Lem. 2.4]{Sh:f}}]
\label{seallem}
Suppose that $T\subseteq {}^{<\omega_1}\omega_1$ is a tree of
height $\omega_1$, $|T|=\aleph_1$, $|\lim(T)|\leq\aleph_1$ and $\bbP_T$ is
the forcing notion for sealing the branches of $T$. 
\begin{enumerate}
\item[(a)] $\bbP_T$ satisfies the ccc.
\item[(b)] If $G\subseteq\bbP_T$ is generic over $\bV$ and $\bV^*$ is a
  universe of {\rm ZFC} extending $\bV[G]$ and such that $(\aleph_1)^{\bV^*} 
  =\aleph_1^\bV (=(\aleph_1)^{\bV[G]})$, then 
\[\bV^*\models\lim(T)=\big(\lim(T)\big)^\bV.\]  
\end{enumerate}
\end{lemma}
\medskip

\noindent{\bf Adding a bound to $\cG\subseteq\bqo$ and a family
  $\cU\subseteq \cP(\omega_1)$:}\quad After sealing branches of a tree, we
  will force a new member $r^*$ of our family $H\subseteq\bqo$ at the same
  time making sure that some family $\cU$ of subsets of $\omega_1$ is
  included in $\fil(r^*)$.   

\begin{definition}
\label{addbddef}
Suppose that $\cG\subseteq\bqo$ and $\cU\subseteq \cP(\omega_1)$ are such
that 
\begin{enumerate}
\item[(a)] $\cG\subseteq\bqo$ is $\leq^*$--directed and 
\item[(b)] $U_0\cap\ldots\cap U_n\in\big(\fil(\cG)\big)^+$ for every
  $U_0,\ldots U_n\in\cU$, $n<\omega$.  
\end{enumerate}
We define a forcing notion $\bbQ^{\rm bd}(\cG,\cU)$ as follows:\\
{\bf a condition $p$ in $\bbQ^{\rm bd}(\cG,\cU)$} is a triple
$(r^p,\bG^p,\bU^p)$ such that $r^p\subseteq\cF^{\rm ult}_{\omega_1}$ is
countable and strongly disjoint (i.e., it satisfies the demands of
\ref{genfil}(2)), $\bG^p\subseteq \cG$ is countable and $\bU^p\subseteq \cU$  
is countable;\\ 
{\bf the order $\leq=\leq_{\bbQ^{\rm bd}(\cG,\cU)}$} is defined by: $p\leq q$
if and only if ($p,q\in \bbQ^{\rm bd}(\cG,\cU)$ and)\\
$\bU^p\subseteq \bU^q$, $\bG^p\subseteq \bG^q$, $r^p\subseteq r^q$ and
for every $(\alpha,Z,d)\in r^q\setminus r^p$ we have that:
\begin{itemize}
\item $\big(\forall (\alpha',Z',d')\in r^p\big)\big(Z'\subseteq\alpha\big)$ and 
\item $\big(\forall r\in\bG^p\big)\big((\alpha,Z,d)\in \Sigma(r)\big)$
  ($\Sigma(r)$ was defined in Definition \ref{bigA}) and 
\item $\big(\forall U\in \bU^p\big)\big(U\cap Z\in d\big)$.
\end{itemize}
\smallskip

\noindent We also define a $\bbQ^{\rm bd}(\cG,\cU)$--name $\name{r}$ by
\qquad 
$\forces_{\bbQ^{\rm bd}(\cG,\cU)}\mbox{`` }\name{r}=\bigcup\big\{r^p:p\in
\Gamma_{\bbQ^{\rm bd}(\cG,\cU)}\big\}\mbox{ ''}$.
\end{definition}

\begin{lemma}
\label{addbdlem}
Assume $\cG\subseteq\bqo$, $\cU\subseteq\cP(\omega_1)$ satisfy demands
(a),(b) of \ref{addbddef}. Then 
\begin{enumerate}
\item $\bbQ^{\rm bd}(\cG,\cU)$ is a $\sigma$--closed forcing notion,
\item $\forces_{\bbQ^{\rm bd}(\cG,\cU)}$`` $\name{r}\in \bqo$ and $(\forall
  r\in\cG)(r\leq^* \name{r})$ and $\cU\subseteq\fil(\name{r})$ ''.  
\end{enumerate}
\end{lemma}

\begin{proof}
(1)\quad Straightforward.
\medskip

\noindent (2)\quad To argue that $\forces_{\bbQ^{\rm bd}(\cG,\cU)}$``
$\name{r}\in \bqo$ '', suppose $p\in \bbQ^{\rm bd}(\cG,\cU)$. Let
$\{r_n:n<\omega\}=\bG^p$, $\{U_n:n<\omega\}=\bU^p$ (we allow
repetitions). Choose inductively $(\alpha_m,Z_m,d_m)\in \cF^{\rm
  ult}_{\omega_1}$ such that for $m<\omega$ we have 
\begin{itemize}
\item $\big(\forall (\alpha',Z',d')\in
  r^p\big)\big(Z'\subseteq\alpha_0\big)$, $Z_m\subseteq\alpha_{m+1}$, and 
\item $(\alpha_m,Z_m,d_m)\in \Sigma(r_0)\cap\ldots\cap\Sigma(r_m)$, and 
\item $U_0\cap\ldots \cap U_m\cap Z_m\in d_m$. 
\end{itemize}
[Why is the choice possible? Since $\cG$ is directed, we may first choose
$s\in\cG$ such that $r_0,\ldots,r_m\leq^* s$. Then for some
$\beta<\omega_1$, if $(\alpha,Z,d)\in s$ and $\beta\leq\alpha$, then
$(\alpha,Z,d)\in \Sigma(r_0)\cap\ldots\cap\Sigma(r_m)$. By the assumption
\ref{addbddef}(b) on $\cU$ we know that $U_0\cap\ldots\cap U_m\cap Z\in d$
for $\omega_1$ many $(\alpha,Z,d)\in s$, so we may choose
$(\alpha_m,Z_m,d_m)\in s$ as required.] 

After the above construction is carried out, pick any uniform ultrafilter
$e$ on $\omega$ and put 
\[\alpha=\alpha_0,\quad Z=\bigcup_{m<\omega} Z_m,\quad \mbox{and}
\quad d=\bigoplus^e_{m<\omega} d_m.\]
Then $q=(r^p\cup\{(\alpha,Z,d)\},\bG^p,\bU^p)\in \bbQ^{\rm bd}(\cG,\cU)$ is
a condition stronger than $p$. Thus by an easy density argument we see that
$\forces_{\bbQ^{\rm bd}(\cG,\cU)}$`` $|\name{r}|=\omega_1$ ''. The rest
should be clear.  
\end{proof}

Let us recall that a very reasonable ultrafilter on $\lambda$ is a weakly
reasonable ultrafilter $D$ such that $D=\fil(H)$ for some
$(<\lambda^+)$--directed family $H\subseteq\bqz$ (see \cite[Def
2.5(5)]{Sh:830}). Now we may state and prove our result. 

\begin{theorem}
  \label{inacneeded}
Assume that $\kappa$ is a strongly  inaccessible cardinal. Then there is a
$\kappa$--cc proper forcing notion $\bbP$ such that   
\[\begin{array}{ll}
\forces_\bbP&\mbox{`` there is a $\leq^*$--increasing sequence $\langle
  r_\xi:\xi<\omega_2\rangle\subseteq\bqo$ such that }\\
&\ \ \fil\big(\{r_\xi:\xi<\omega_2\}\big)\mbox{ is a very reasonable
  ultrafilter on }\omega_1\\
&\ \mbox{ but Odd has a winning strategy in the game
}\Game_{\{r_\xi:\xi<\omega_2\}}\mbox{ ''.}
\end{array}\]    
\end{theorem}

\begin{proof}
The forcing notion $\bbP$ will be obtained as the limit of a CS iteration
of proper forcing notions $\langle\bbP_\xi,\name{\bbQ}_\xi:\xi<\kappa
\rangle$.  The iteration will be built so that for each $\xi<\kappa$ 
\[\forces_{\bbP_\xi}\mbox{`` } \name{\bbQ}_\xi\mbox{ is a proper
  forcing notion of size $<\kappa$ ''},\]
so we will be sure that the intermediate stages $\bbP_\xi$ and the limit
$\bbP_\kappa$ will be proper and each $\bbP_\xi$ (for $\xi<\kappa$) will
have a dense subset of cardinality $<\kappa$. Thus $\bbP_\kappa$ will
satisfy $\kappa$--cc (and $\kappa$ will not be collapsed). Since in the
process of iteration we will also collapse to $\aleph_1$ all uncountable
cardinals below $\kappa$,  we will know that 
\[\forces_{\bbP_\kappa}\mbox{`` }\aleph_1=(\aleph_1)^\bV\ \&\ 2^{\aleph_1}=
\aleph_2=\kappa\mbox{ ''}.\] 
Thus we may set up a bookkeeping device that gives us a list $\langle
\name{C}_\zeta, \name{A}_\zeta,\name{\rho}_\zeta:\zeta<\kappa\rangle$ such
that 
\begin{itemize}
\item $\name{C}_\zeta$ is a $\bbP_\zeta$--name for a club of $\omega_1$,
\item $\name{A}_\zeta$ is a $\bbP_\zeta$--name for a subset of $\omega_1$,
\item $\name{\rho}_\zeta$ is a $\bbP_\zeta$--name for a function from
  $\omega_1$ to   $\omega_1$,  and 
\item for each $\bbP_\kappa$--name $\name{C}$ for a club of $\omega_1$, for
some $\zeta<\kappa$ we have $\forces_{\bbP_\kappa}\name{C}=\name{C}_\zeta$,
and similarly for names $\name{A}$ for subsets of $\omega_1$ and names
$\name{\rho}$ for elements of ${}^{\omega_1}\omega_1$. 
\end{itemize}
\medskip

Before continuing let us set some terminology used later. {\em A partial
  strategy\/} is a function $\sigma$  such that  
\begin{itemize}
\item $\dom(\sigma)\subseteq\{\eta\in{}^{<\omega_1}\omega_1:\lh(\eta)$ is 
  an odd ordinal $\}$, and
\item $\big(\forall\nu\in\dom(\sigma)\big)\big(\sigma(\nu) \in\omega_1
  \setminus(\sup(\nu)+1)\big)$.
\end{itemize}
We say that a sequence $\eta\in{}^{\leq\omega_1}\omega_1$ {\em is played
  according to a partial strategy $\sigma$} if    
\begin{itemize}
\item the sequence $\eta$ is increasing continuous, and 
\item for every odd ordinal $\alpha<\lh(\eta)$ we have $\eta\rest\alpha \in
  \dom(\sigma)$ and $\eta(\alpha)=\sigma(\eta\rest\alpha)$. 
\end{itemize}
If $\rho,\eta\in {}^{\omega_1}\omega_1$ and $\eta$ is played according to
$\sigma$, then we say that $\eta=\sigma[\rho]$ if $\eta(0)=\rho(0)$ and
$\eta(2\alpha+2)=\eta(2\alpha+1)+\rho(1+\alpha)+1$ for each
$\alpha<\lambda$. Also, for an increasing sequence $\eta\in{}^{\omega_1}
\omega_1$ let   
\[U_\eta=\bigcup\big\{\big[\eta(2\alpha),\eta(2\alpha+1)\big): \alpha<
\omega_1\big\}.\]  
\medskip

Now, we will inductively choose $\name{\bbQ}_\xi$ and $\name{T}_\xi,
\name{\sigma}_\xi,\name{r}_\xi$ so that for each $\xi<\kappa$ the following
demands are satisfied.  
\begin{enumerate}
\item[$(\boxplus)_1$]  $\name{r}_\xi$ is a $\bbP_{\xi+1}$--name for a member
  of $\bqo$ and $\forces_{\bbP_{\xi+1}}(\forall\zeta<\xi)(\name{r}_\zeta
  \leq^*\name{r}_\xi)$,
\item[$(\boxplus)_2$] $\name{T}_\xi$ is a $\bbP_\xi$--name for a subtree of
  ${}^{<\omega_1}\omega_1$ of height $\omega_1$ (with no maximal nodes). 
\item[$(\boxplus)_3$] $\name{\sigma}_\xi$ is a $\bbP_\xi$--name for a
  partial strategy with domain $\{\eta\in\name{T}_\xi:\lh(\eta)\mbox{ is odd
  }\}$, and all nodes of the tree $\name{T}_\xi$ are played according to 
  $\name{\sigma}_\xi$.
\item[$(\boxplus)_4$] $\forces_{\bbP_{\xi+1}}\big(\exists\eta\in\lim(
  \name{T}_{\xi+1})\big)\big(\eta=\name{\sigma}_{\xi+1}[
  \name{\rho}_\xi]\big)$. 
\item[$(\boxplus)_5$] $\forces_{\bbP_\xi}(\forall \zeta<\xi)(\name{T}_\zeta
  \subseteq \name{T}_\xi\ \&\ \name{\sigma}_\zeta\subseteq
  \name{\sigma}_\xi)$ and 
\[\begin{array}{ll}
\forces_{\bbP_{\xi+1}}&\mbox{`` if }\nu\in {}^{<\omega_1}\omega_1\mbox{ is  
increasing continuous and such that}\\
&\ \ \lh(\nu)=\gamma+1\mbox{ for a limit $\gamma$ and }(\forall\alpha<
\gamma)(\nu\rest\alpha\in\name{T}_\xi)\mbox{ but }\nu\rest\gamma\notin
\name{T}_\xi,\\  
&\ \mbox{ then }\nu\in \name{T}_{\xi+1}\mbox{ ''.} 
  \end{array}\] 
\item[$(\boxplus)_6$] $\forces_{\bbP_\xi}\big(\forall\eta_0,\ldots,\eta_n
  \in\lim(\name{T}_\xi)\big)\big(\forall\zeta<\xi\big)\big(U_{\eta_0}\cap
  \ldots\cap U_{\eta_n}\in\big(\fil(\name{r}_\zeta)\big)^+\big)$ for each
  $n<\omega$, and\\
 $\forces_{\bbP_{\xi+1}}\big(\forall\eta\in\lim(
\name{T}_\xi)\big)\big(U_\eta\in\fil(\name{r}_\xi)\big)$.
\item[$(\boxplus)_7$] $\forces_{\bbP_{\xi+1}}$`` $\name{A}_\xi\in
  \fil(\name{r}_\xi)$ or $\omega_1\setminus\name{A}_\xi\in\fil(
  \name{r}_\xi)$ '' and $\forces_{\bbP_{\xi+1}}$`` if $\langle\delta_\alpha:
  \alpha<\lambda\rangle$ is the increasing enumeration of
  $\name{C}_\xi$, then for some club $C^*\subseteq\omega_1$ we have
  $\omega_1\setminus\bigcup\big\{[\delta_\alpha,\delta_{\alpha+1}):\alpha
  \in C^*\big\}\in \fil(\name{r}_\xi)$ ''.  
\item[$(\boxplus)_8$] For $\xi>0$, $\name{\bbQ}_\xi$ is the $\bbP_\xi$--name
  for the composition 
\[\bbC*\name{\bbL}*\name{\bbP}_{\name{T}_\xi}*\name{\bbQ}^{\rm 
bd}(\{\name{r}_\zeta:\zeta<\xi\},\{U_\eta:\eta\in \lim(\name{T}_\xi)\})\] 
(see \ref{CLpres}, \ref{sealdef}, \ref{addbddef}). Hence we know that also 
\item[$(\boxplus)_9$] for every $\bbP_{\xi+1}$--name $\name{\bbQ}$ for a
  proper forcing notion,
$\forces_{\bbP_{\xi+1}*\name{\bbQ}}\mbox{`` }\lim(\name{T}_\xi)=
\big(\lim(\name{T}_\xi)\big)^{\bV^{\bbP_\xi}}\mbox{ ''}$.
\end{enumerate}
\medskip

To start, we let $r_{-1}$ be any fixed element of $\bqo$. We choose 
$\sigma':{}^{<\omega_1}\omega_1\longrightarrow\omega_1$ so that for every 
$\eta\in {}^{<\omega_1}\omega_1$ there is $(\alpha,Z,d)\in r$ such that
$\sup(\eta)<\alpha$ and $Z\subseteq \sigma'(\eta)$, and we let
$\name{T}_0=T_0=\{\sigma'[\rho_0]\rest\alpha:\alpha<\omega_1\} \subseteq 
{}^{<\omega_1}\omega_1$.  (So $T_0$ is a tree with
$\lim(T_0)=\{\sigma'[\rho_0]\}$.) Finally $\name{\sigma}_0=\sigma_0=
\sigma'\rest \{\nu\in T_0:\lh(\nu)\mbox{ is odd }\}$.  Now, the forcing
notion $\bbQ_0$ is:  
\[\bbC*\name{\bbL}*\name{\bbP}_{T_0}*\name{\bbQ}^{\rm bd}(\{r_{-1}\},
\{U_{\sigma_0[\rho_0]}\}).\]   
Clearly, the families $\{r_{-1}\}$ and $\{U_{\sigma_0[\rho_0]}\}$ satisfy
the demands (a),(b) of Definition \ref{addbddef}.   
\medskip

Now suppose that we have arrived to a successor stage $\xi=\zeta+1$ (and we
have already defined $\bbP_\zeta$ and $\bbP_\zeta$--names $\name{T}_\zeta,
\name{\sigma}_\zeta$, and $\bbP_{\vare+1}$--names $\name{r}_\vare$ for
$\vare<\zeta$ so that the demands of $(\boxplus)_1$--$(\boxplus)_6$ hold. It
follows from $(\boxplus)_1+(\boxplus)_6$ that $\name{\bbQ}_\zeta$ is
correctly determined by clause $(\boxplus)_8$, so   
\[\forces_{\bbP_\zeta}\name{\bbQ}_\zeta= \bbC*\name{\bbL}* 
\name{\bbP}_{T_\zeta} *\name{\bbQ}^{\rm bd}\big(\{\name{r}_\vare: \vare<
\zeta\}, \{U_\eta: \eta\in\lim(\name{T}_\zeta)\}\big).\]  
(Remember also that, by $(\boxplus)_9$, all $\omega_1$--branches of
$\name{T}_\zeta$ in extensions by proper forcing over
$\bV^{\bbP_\zeta*\name{\bbQ}_\zeta}$ are the same as those in
$\bV^{\bbP_\zeta}$.) Note, that the last factor of $\name{\bbQ}_\zeta$ adds
an element $\name{r}\in \bqo$ (see \ref{addbdlem}(2)) and we know that 
\[\forces_{\bbP_\zeta*\name{\bbQ}_\zeta}\mbox{`` }\big(\forall\vare<\zeta
\big)\big(\name{r}_\vare\leq^*\name{r}\big)\mbox{ and }
\{U_\eta: \eta\in\lim(\name{T}_\zeta)\}\subseteq\fil(\name{r})\mbox{
  ''.}\]  
In $\bV^{\bbP_\zeta*\name{\bbQ}_\zeta}$, we may choose thin enough
uncountable subset of $\name{r}$, getting $\name{r}'\subseteq\name{r}$
satisfying the demand in $(\boxplus)_7$ and such that 
\[\big(\forall (\alpha,Z,d),(\alpha',Z',d') \in\name{r}'\big)\big(\alpha<
\alpha'\ \Rightarrow\ \sup(Z)+\omega<\alpha'\big).\]
Let $\name{\sigma}':{}^{<\omega_1}\omega_1\longrightarrow\omega_1$ be such
that $\name{\sigma}'\rest\dom(\name{\sigma}_\zeta)=\name{\sigma}_\zeta$ and
for $\nu\in {}^{<\omega_1}\omega_1\setminus \dom(\name{\sigma}_\zeta)$ we
have 
\begin{enumerate}
\item[$(\circledast)_1$] \quad $\name{\sigma}'(\nu)=\min\big\{\beta
<\omega_1:\big(\exists (\alpha,Z,d) \in\name{r}'\big)\big(\sup(\nu)<\alpha\
\&\ Z\subseteq\beta\big)\big\}$. 
\end{enumerate}
Let $\name{\eta}^*=\name{\sigma}'[\name{\rho}_\zeta]$ and let
$\name{r}_\zeta=\big\{(\alpha,Z,d)\in\name{r}':U_{\name{\eta}^*}\cap Z\in 
d\big\}$. It follows from our choices so far that $\name{r}_\zeta\in\bqo$,
and $\name{r}_\vare\leq^*\name{r}_\zeta$ for $\vare<\zeta$, and also  
\begin{enumerate}
\item[$(\circledast)_2$] \quad for each $i<\omega_1$, 
\[\name{\eta}^*\rest (2i+1)\notin \name{T}_\zeta\ \Rightarrow\ \big(\exists 
  (\alpha,Z,d) \in\name{r}_\zeta\big)\big(\name{\eta}^*(2i)<\alpha\ \&\ Z
  \subseteq\name{\eta}^*(2i+1)\big).\]
\end{enumerate}
Put $\name{T}^*_\zeta=\name{T}_\zeta\cup \{\name{\eta}^*\rest \alpha:
\alpha<\omega_1\}$ and define $\name{\sigma}'':{}^{<\omega_1}\omega_1
\longrightarrow\omega_1$ so that $\name{\sigma}''\rest\name{T}^*_\zeta=
\name{\sigma}'\rest \name{T}^*_\zeta$ and for $\nu\in
{}^{<\omega_1}\omega_1\setminus \name{T}^*_\zeta$ we have
\begin{enumerate}
\item[$(\circledast)_3$] \quad $\name{\sigma}''(\nu)=\min\big\{\beta
<\omega_1:\big(\exists (\alpha,Z,d)\in\name{r}_\zeta\big)\big(\sup(\nu)
<\alpha\ \&\ Z\subseteq\beta\big)\big\}$.  
\end{enumerate}
Let 
\[\begin{array}{ll}
\name{S}=\{\nu\in{}^{<\omega_1}\omega_1:&\nu\mbox{ is increasing continuous
  of length }\lh(\nu)=\gamma+1\\
&\mbox{for some limit }\gamma\mbox{ and }(\forall\alpha<\gamma)(\nu\rest
\alpha\in \name{T}_\zeta)\mbox{ but }\nu\rest\gamma\notin\name{T}^*_\zeta
\}.  
  \end{array}\]
For each $\nu\in\name{S}$ let $\eta_\nu\in {}^{\omega_1}\omega_1$ be such
that $\nu\vtl\eta_\nu$, $\eta_\nu$ is played according to $\name{\sigma}''$
and for every odd ordinal $\alpha\geq\lh(\nu)$ we have $\eta_\nu(\alpha+1)= 
\eta_\nu(\alpha)+889$.  Put  $\name{T}_{\zeta+1}= \name{T}_\zeta^*\cup
\{\eta_\nu\rest\alpha:\nu\in \name{S}\ \&\ \alpha<\omega_1\}$.
Note that (still in $\bV^{\bbP_\zeta*\name{\bbQ}_\zeta}$) we have that
$\lim\big(\name{T}_{\zeta+1}\big)=\lim\big(\name{T}_\zeta\big)\cup 
\{\eta_\nu:\nu\in\name{S}\} \cup\{\name{\eta}^*\}$.

It follows from the choice of $\name{r}$, $\name{r}_\zeta$ that
$\big(\forall\eta\in\lim(\name{T}_\zeta)\big)\big(U_\eta\in \fil(
\name{r}_\zeta)\big)$ and by the definition of $\name{\sigma}''$ we get that
$\big(\forall\nu \in\name{S}\big)\big(U_{\eta_\nu}\in\fil(\name{r}_\zeta)
\big)$ (remember the choice of $\name{r}'$). Hence, remembering the 
definition of $\name{r}_\zeta$, we conclude that $\big(\forall\eta\in
\lim(T_{\zeta+1})\big)\big(U_\eta\in \fil(\name{r}_\zeta) \big)$.

Finally we put $\name{\sigma}_{\zeta+1}=\name{\sigma}''\rest \{\nu\in 
T_{\zeta+1}: \lh(\nu)\mbox{ is odd }\}$. One easily verifies that the
relevant demands in $(\boxplus)_1$--$(\boxplus)_7$ hold for
$\name{T}_{\zeta+1}$, $\name{\sigma}_{\zeta+1}$ and $\name{r}_\zeta$. Let us
also stress for future reference that 
\begin{enumerate}
\item[$(\circledast)_4$] if $\nu\in\dom(\name{\sigma}_{\zeta+1})\setminus
  \name{T}_\zeta$ is of length $2i+1$, then there is $(\alpha,Z,d)\in
  \name{r}_\zeta$ such that $\nu(2i)<\alpha$ and $Z\subseteq
  \name{\sigma}_{\zeta+1}(\nu)$.  
\end{enumerate}
\medskip

Suppose now that we have arrived to a limit stage $\xi<\kappa$ and we have
defined $\bbP_\zeta$ names $\name{\bbQ}_\zeta,\name{T}_\zeta,
\name{\sigma}_\zeta$, and $\name{r}_\zeta$ for all $\zeta<\xi$. In
$\bV^{\bbP_\xi}$ we define $\name{T}_\xi=\bigcup\limits_{\zeta<\xi}
\name{T}_\zeta$ and $\name{\sigma}_\xi=\bigcup\limits_{\zeta<\xi}
\name{\sigma}_\zeta$. We have to argue that the relevant demands in 
$(\boxplus)_2$--$(\boxplus)_6$ are satisfied, and the only problematic one
is the first condition of $(\boxplus)_6$. If $\cf(\xi)=\omega_0$, then
$\forces_{\bbP_\xi}\lim(\name{T}_\xi)=\bigcup_{\zeta<\xi}\lim(
\name{T}_\zeta)$, so there are no problems. We will show that $(\boxplus)_6$ 
holds also if $\cf(\xi)\geq\omega_1$ and for this we will argue {\em a
  contrario}. 

Suppose towards contradiction that ($\cf(\xi)\geq \omega_1$ and) we have
$\bbP_\xi$--names $\name{\eta}_0,\ldots,\name{\eta}_n$ ($n<\omega$) and a
condition $p\in\bbP_\xi$ such that 
\[p\forces_{\bbP_\xi}\mbox{`` }\name{\eta}_0,\ldots,\name{\eta}_n\in
\lim(\name{T}_\xi) \mbox{ and }\big(\exists\zeta<\xi\big)\big(\omega_1
\setminus (U_{\name{\eta}_0}\cap\ldots \cap U_{\name{\eta}_n})\in
\fil(\name{r}_\zeta)\big)\mbox{ ''.}\] 
Remembering that $(\boxplus)_1+(\boxplus)_6$ hold on earlier stages, we may
pass to a stronger condition (if necessary) and assume additionally that for
some $\zeta<\xi$, $\gamma<\omega_1$ and pairwise distinct
$\nu_0,\ldots,\nu_n\in {}^\gamma\omega_1$ we have 
\[\begin{array}{ll}
p\forces_{\bbP_\xi}&\mbox{`` }\name{\eta}_0,\ldots,\name{\eta}_n\notin 
\bigcup\limits_{\vare<\xi}\lim(\name{T}_\vare)\mbox{ and }\name{\eta}_0
\rest \gamma=\nu_0,\ldots,\name{\eta}_n\rest\gamma=\nu_n\mbox{ and }\\
&\ \ \big(\forall (\alpha,Z,d)\in \name{r}_\zeta \big)\big(\gamma\leq\sup(Z)\
\Rightarrow\ U_{\name{\eta}_0}\cap\ldots \cap U_{\name{\eta}_n}\cap Z\notin d\big)\mbox{ ''.}
\end{array}\] 
The forcing notion $\bbP_\xi$ is proper, so we may choose a countable
elementary submodel $N\prec \cHchi$ such that $\name{\eta}_0, \ldots,
\name{\eta}_n, \nu_0,\ldots,\nu_n,\zeta,\xi, \gamma,p,\ldots \in N$ and then
we may pick an $(N,\bbP_\xi)$--generic condition $q\geq p$. Let
$\gamma^*=N\cap\omega_1$ and $\xi^*=\sup(N\cap \xi)$, and we may assume
$q\in\bbP_{\xi^*}$. Then   
\begin{enumerate}
\item[$(\circledast)_5$] $q\forces_{\bbP_{\xi^*}}$`` $(\forall
  i\leq n)(\forall \vare<\xi^*)(\exists \delta<\gamma^*)(\name{\eta}_i\rest
  \delta\notin \name{T}_\vare)$ '',  and hence 
$q\forces_{\bbP_{\xi^*}}$`` $\name{\eta}_0\rest\gamma^*,\ldots,\name{\eta}_n
\rest\gamma^*\notin \name{T}_{\xi^*}$ ''.  
\end{enumerate}
Why? As  for each $\vare\in N\cap\xi$ we have a name $\name{\delta}\in
N$ for an ordinal below $\omega_1$ such that $p\forces \name{\eta}_i\rest
\name{\delta}\notin \name{T}_\vare$, so we may use the genericity of $q$. By
a similar argument, 
\begin{enumerate}
\item[$(\circledast)_6$] $q\forces_{\bbP_{\xi^*}}$`` $(\forall i\leq n)
  (\forall\delta<\gamma^*)(\exists \vare<\xi^*)(\name{\eta}_i\rest \delta\in
  \name{T}_\vare)$ '', so also
$q\forces_{\bbP_{\xi^*}}$`` $(\forall\delta<\gamma^*)(\name{\eta}_0\rest
\delta,\ldots,\name{\eta}_n\rest\delta\in\name{T}_{\xi^*})$ '',
\end{enumerate}
and
\begin{enumerate}
\item [$(\circledast)_7$] $q\forces_{\bbP_{\xi^*}}$`` $\name{\eta}_0(
  \gamma^*)= \ldots=\name{\eta}_n(\gamma^*)=\gamma^*$ ''.
\end{enumerate}
(Remember that $\name{\eta}_i$ are increasing continuous.) Now, consider a
$\bbP_{\xi^*}$--name $\name{q}$ for the following member of $\bbQ_{\xi^*}$: 
\[\big(\emptyset_{\bbC},\emptyset_{\name{\bbL}},
\emptyset_{\bbP_{\name{T}_{\xi^*}}}, (\emptyset,\{\name{r}_\zeta\},
\emptyset) \big).\]
Directly from the definition of the order of the forcing $\bbQ^{\rm bd}$ and
the choice of $\name{r}_{\xi^*}$ we see that
\[q\cup\{(\xi^*,\name{q})\}\forces_{\bbP_{\xi^*+1}}\mbox{`` } \name{r}_{\xi^*}
\subseteq\Sigma(\name{r}_\zeta)\mbox{ ''.}\] 
It follows from $(\circledast)_5,(\circledast)_6$ and $(\boxplus)_5$ that 
\[q\cup\{(\xi^*,\name{q})\}\forces_{\bbP_{\xi^*+1}}\mbox{`` } \name{\eta}_0
\rest(\gamma^*+1), \ldots,\name{\eta}_n\rest(\gamma^*+1) \in
\name{T}_{\xi^*+1}\setminus\name{T}_{\xi^*}\mbox{ '',}\]   
so let us look what are the respective values of the partial strategy
$\name{\sigma}_{\xi^*+1}$. By $(\circledast)_4$ we know that  
\[\begin{array}{ll}
q\cup\{(\xi^*,\name{q})\}\forces_{\bbP_{\xi^*+1}}&\mbox{`` }\mbox{ there
  exists }(A,Z,d)\in\name{r}_{\xi^*}\mbox{ such that for each } i\leq n\\
&\quad\name{\eta}_i(\gamma^*)=\gamma^*<\alpha\mbox{ and }Z \subseteq 
\name{\eta}_i(\gamma^*+1)\mbox{ ''.}\end{array}\]
Since $\gamma^*>\gamma$ we get a contradiction with the choice of $p$.
\medskip

This completes the inductive definition of the iteration and the names
$\name{T}_\xi$, $\name{\sigma}_\xi$ and $\name{r}_\xi$. It should be clear
that
\[\begin{array}{ll}
\forces_{\bbP_\kappa}&\mbox{`` the sequence }\langle\name{r}_\xi:\xi<\kappa
\rangle\mbox{ is $\leq^*$--increasing and}\\
&\ \ \fil\big(\{r_\xi:\xi<\omega_2\}\big)\mbox{ is a very reasonable
  ultrafilter on }\omega_1\mbox{ and }\\
&\ \ \bigcup\limits_{\xi<\kappa}\name{\sigma}_\xi\mbox{ is a winning
  strategy for Odd in the game }\Game_{\{r_\xi:\xi<\kappa\}}\mbox{ ''.}
\end{array}\]    
\end{proof}


\begin{thebibliography}{1}

\bibitem{J}
Thomas Jech.
\newblock {\em {Set theory}}.
\newblock Springer Monographs in Mathematics. Springer-Verlag, Berlin, 2003.
\newblock The third millennium edition, revised and expanded.

\bibitem{MRSh:799}
Pierre Matet, Andrzej Roslanowski, and Saharon Shelah.
\newblock {Cofinality of the nonstationary ideal}.
\newblock {\em Transactions of the American Mathematical Society},
  357:4813--4837, 2005.
\newblock math.LO/0210087.

\bibitem{RoSh:890}
Andrzej Roslanowski and Saharon Shelah.
\newblock {Reasonable ultrafilters, again}.
\newblock {\em Notre Dame Journal of Formal Logic}, submitted.
\newblock math.LO/0605067.

\bibitem{RoSh:470}
Andrzej Roslanowski and Saharon Shelah.
\newblock {Norms on possibilities I: forcing with trees and creatures}.
\newblock {\em {Memoirs of the American Mathematical Society}}, 141(671):xii +
  167, 1999.
\newblock math.LO/9807172.

\bibitem{Sh:288}
Saharon Shelah.
\newblock {Strong Partition Relations Below the Power Set: Consistency, Was
  Sierpi\'nski Right, II?}
\newblock In {\em {Proceedings of the Conference on Set Theory and its
  Applications in honor of A.Hajnal and V.T.Sos, Budapest, 1/91}}, volume~60 of
  {\em Colloquia Mathematica Societatis Janos Bolyai. Sets, Graphs, and
  Numbers}, pages 637--638. 1991.
\newblock math.LO/9201244.

\bibitem{Sh:f}
Saharon Shelah.
\newblock {\em {Proper and improper forcing}}.
\newblock {Perspectives in Mathematical Logic}. {Springer}, 1998.

\bibitem{Sh:546}
Saharon Shelah.
\newblock {Was Sierpi\'nski right? IV}.
\newblock {\em {Journal of Symbolic Logic}}, 65:1031--1054, 2000.
\newblock math.LO/9712282.

\bibitem{Sh:830}
Saharon Shelah.
\newblock {The combinatorics of reasonable ultrafilters}.
\newblock {\em Fundamenta Mathematicae}, 192:1--23, 2006.
\newblock math.LO/0407498.

\end{thebibliography}
\end{document}